\documentclass[showkeys,floatfix,aps,nofootinbib,10pt,a4paper]{revtex4-1}
\usepackage{natbib}
\usepackage{amssymb}
\usepackage{amsmath,amsfonts}
\usepackage{setspace}

\usepackage{listings}
\usepackage{multirow}
\usepackage{epsfig}
\usepackage{graphicx}
\usepackage{mathtools}
\usepackage[caption=false]{subfig}
\graphicspath{{figures/}}
\usepackage{dcolumn}
\begin{document}
\title{Computing ultra-precise eigenvalues of the Laplacian within polygons}
\author{\firstname{Robert} Stephen \surname{Jones}}
\email{\url{mailto:rsjones7@yahoo.com}}
\homepage{\url{http://www.hbeLabs.com}}
\affiliation{Independent Researcher; Changsha, Hunan, China}

\keywords{Laplacian eigenvalue; Helmholtz equation; Method of
  Particular Solutions; point-matching method; polygon; eigenvalue
  bound}

\begin{abstract}
The main difficulty in solving the Helmholtz equation within polygons
is due to non-analytic vertices. By using a method nearly identical to
that used by Fox, Henrici, and Moler in their 1967 paper; it is
demonstrated that such eigenvalue calculations can be extended to
unprecedented precision, very often to well over a hundred digits, and
sometimes to over a thousand digits.

A curious observation is that as one increases the number of terms in
the eigenfunction expansion, the approximate eigenvalue may be made to
alternate above and below the exact eigenvalue. This alternation
provides a new method to bound eigenvalues, by inspection.

Symmetry must be exploited to simplify the geometry, reduce the number
of non-analytic vertices and disentangle degeneracies.  The
symmetry-reduced polygons considered here have at most one
non-analytic vertex from which all edges can be seen.  Dirichlet,
Neumann, and periodic-type edge conditions, are independently imposed
on each polygon edge.

The full shapes include the regular polygons and some
with re-entrant angles (cut-square, L-shape, 5-point
star). Thousand-digit results are obtained for the lowest Dirichlet
eigenvalue of the L-shape, and regular pentagon and hexagon.
\end{abstract}
\maketitle
\section*{\label{sec:intro}Introduction}
\begin{figure}[htb!]
\includegraphics[scale=1.2]{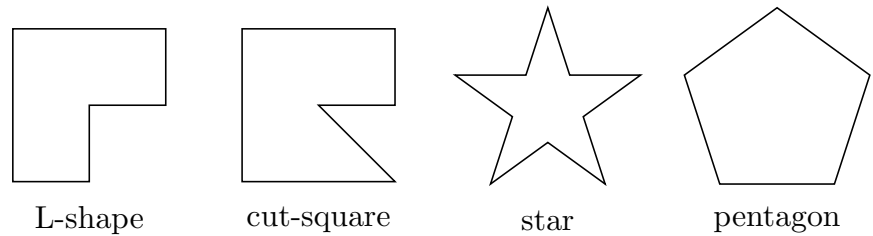}
\caption{Four geometries (with abbreviated names) in which a variety
  of Dirichlet and Neumann eigenvalues are calculated to within at
  least 100 digits. Of the regular polygons, only the regular pentagon
  is shown.}
\label{fig:allshapes}
\end{figure}

The task is to calculate very precise eigenvalues of the Laplacian
within the shapes shown in Fig.~\ref{fig:allshapes}, on which may
be imposed either Neumann or Dirichlet boundary conditions.

The technique is substantially identical to the method used by Fox,
Henrici, and Moler~\cite{fhm1967} (hereafter referred to as ``FHM'')
who used a method of particular solutions (MPS) called the
``point-matching'' or ``collocation'' method to calculate Dirichlet
eigenvalues within the now-famous L-shape. 

To see what is possible using this method, an assorted set of
eigenvalues, all truncated to 100 digits of precision, for the chosen
shapes (Fig.~\ref{fig:allshapes}) is presented in
Table~\ref{tab:ultraprecise}.  In addition, three ``thousand-digit''
results are submitted to the ``On-line Encyclopedia of Integer
Sequences''~\cite{oeisorg} (OEIS.org).  These results far exceed all
previous published results.\footnote{Very recently, P.~Amore,
et~al.~\cite{abfr2015} have independently calculated several highly
precise eigenvalues (up to dozens of digits) for various shapes,
including the L-shape and the cut-square. They also used the FHM
method for some problems, but their focus was on a high-order
Richardson extrapolation method with finite elements, which, among
other things, can handle more varied shapes. Fortunately, their
calculations provide (partial) independent validation of my results.}

\begin{figure}[htb!]
\subfloat[NB-1 (even and odd) $\lambda\approx 8.143$]{\includegraphics[scale=0.7]{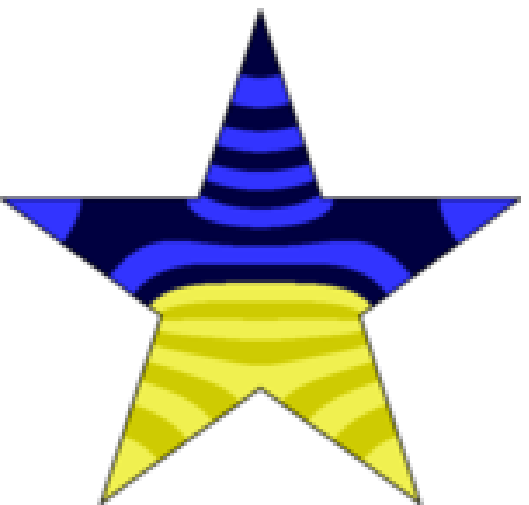}~\includegraphics[scale=0.7]{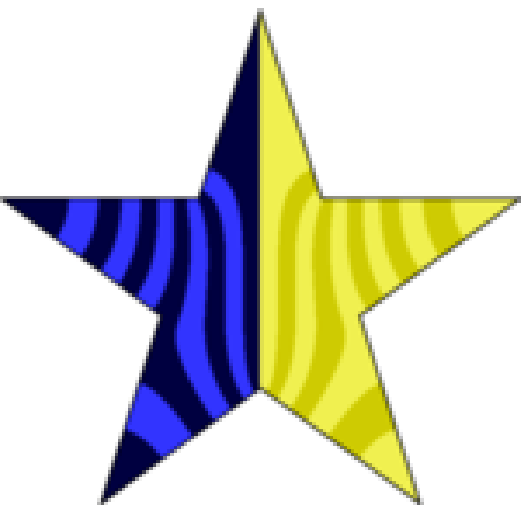}}%
\qquad\qquad\subfloat[NC-1 (even and odd) $\lambda\approx 12.399$]{\includegraphics[scale=0.7]{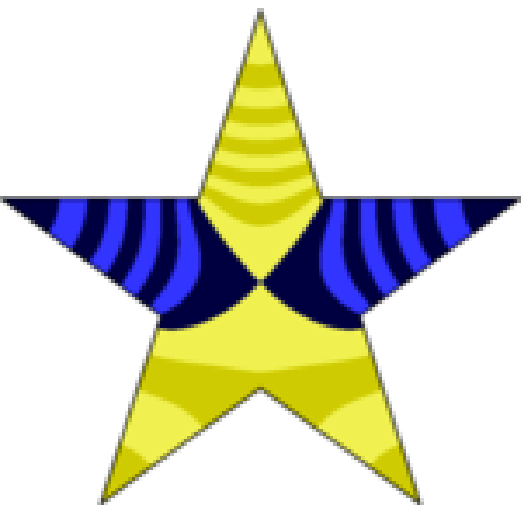}~\includegraphics[scale=0.7]{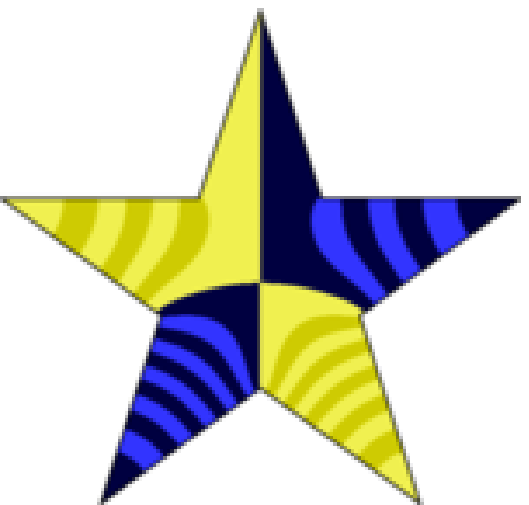}}

\subfloat[NB-12 (even and odd) $\lambda\approx 647.5$]{\includegraphics[scale=0.7]{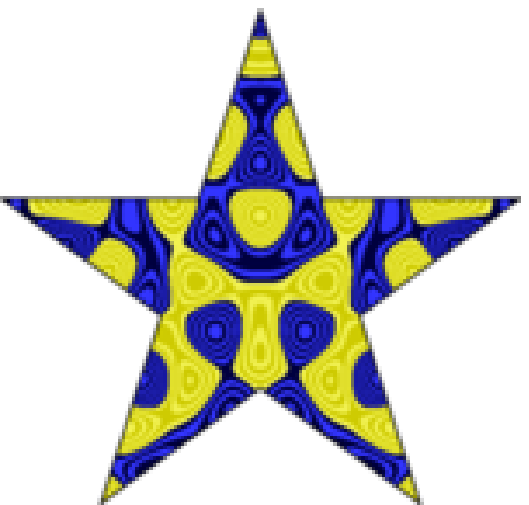}~\includegraphics[scale=0.7]{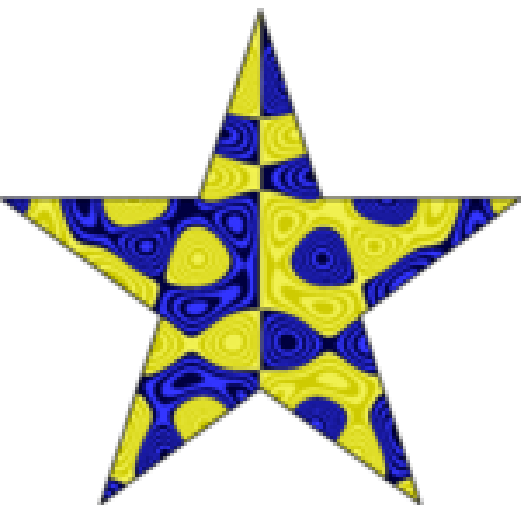}}%
\qquad\qquad\subfloat[DS-18 $\lambda\approx 1073.6$ and DS-9 $\lambda\approx 1080.9$\label{fig:starandpentagon}]{\includegraphics[scale=0.7]{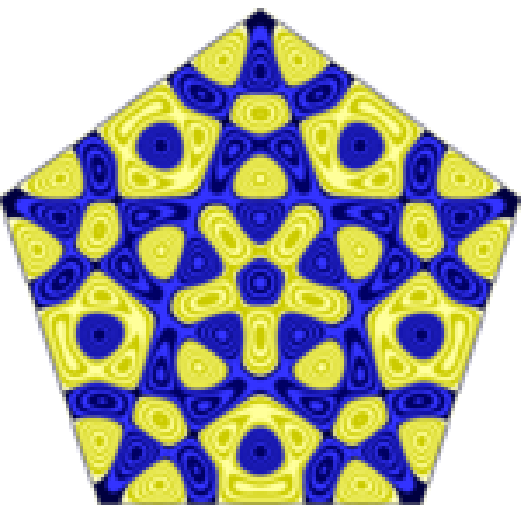}~\includegraphics[scale=0.7]{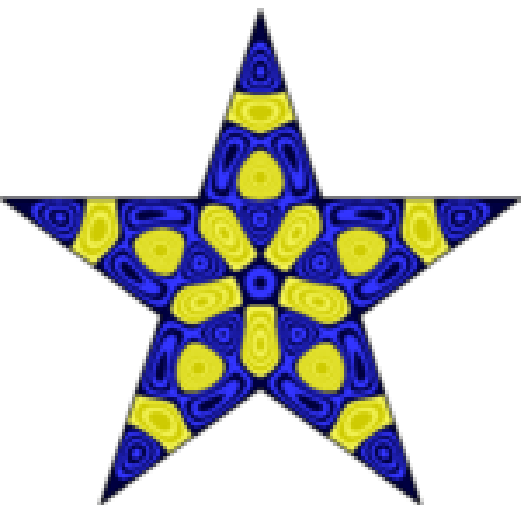}}
\caption{Example contour plots.  Blue and yellow areas are opposite
  sign, and light-dark steps indicate level
  contours. \mbox{(a)--(c)~Three} sets of degenerate Neumann star
  pairs. \mbox{(d)~(LEFT) A} relatively low regular pentagon Dirichlet
  eigenfunction clearly revealing a {\em pentagram\/} nodal
  pattern, and (RIGHT) the corresponding star Dirichlet eigenfunction.}
\label{fig:stareigenfunctions}
\end{figure}

\begin{figure}[!htb]
  \subfloat[\mbox{DC-27 (even and odd) $\lambda\approx 1044.72$}\label{fig:DS027a}]{\includegraphics[scale=0.52]{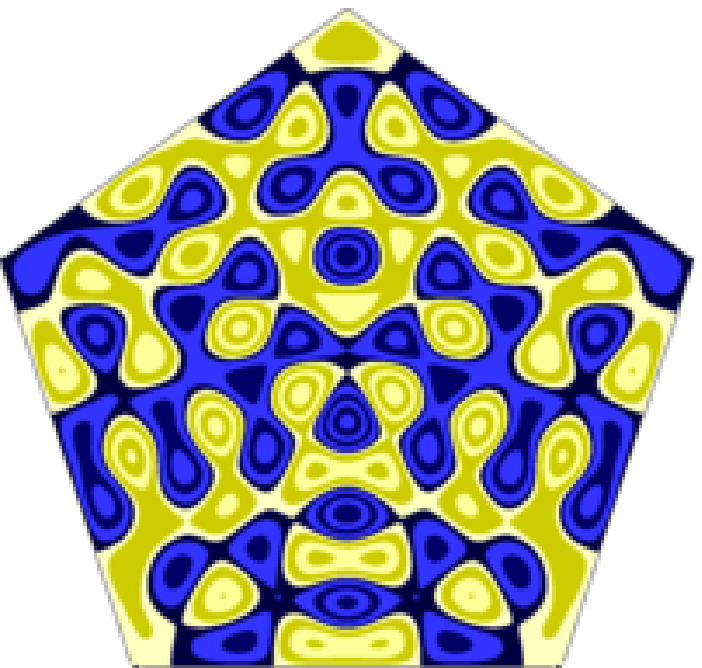}\quad%
    \includegraphics[scale=0.52]{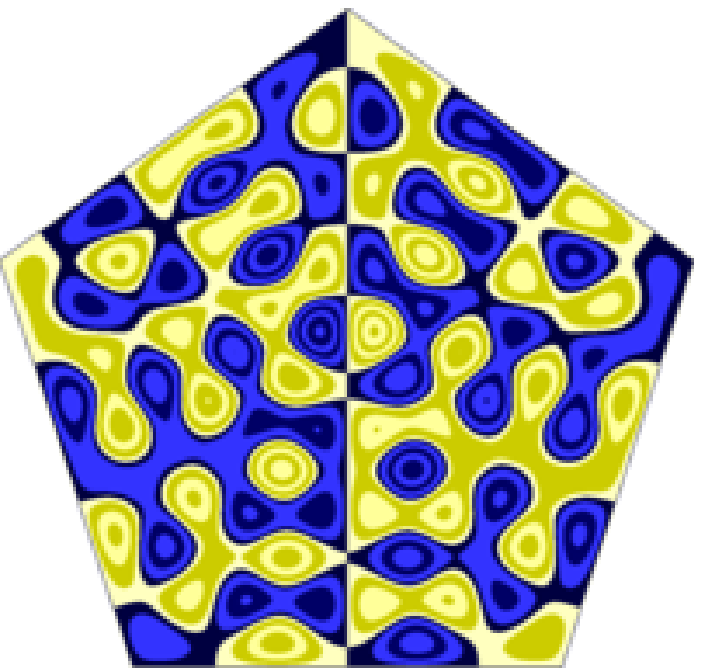}}\qquad%
  \subfloat[\mbox{DS-65: $\lambda\approx 4293.91$}\label{fig:DS065a}]{\includegraphics[scale=0.7]{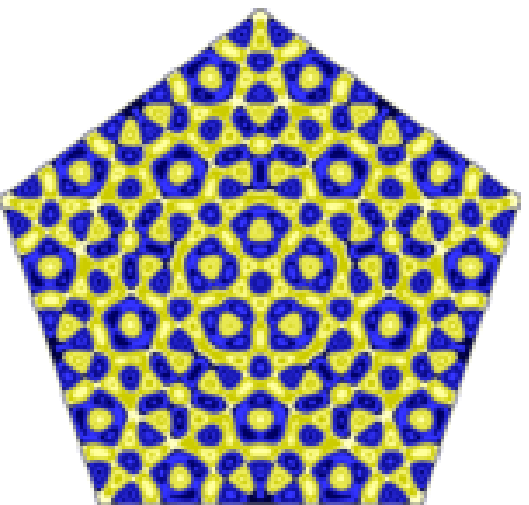}}\qquad%
\subfloat[\mbox{DS-1396: $\lambda\approx 100001.28$}]{\includegraphics[scale=0.52]{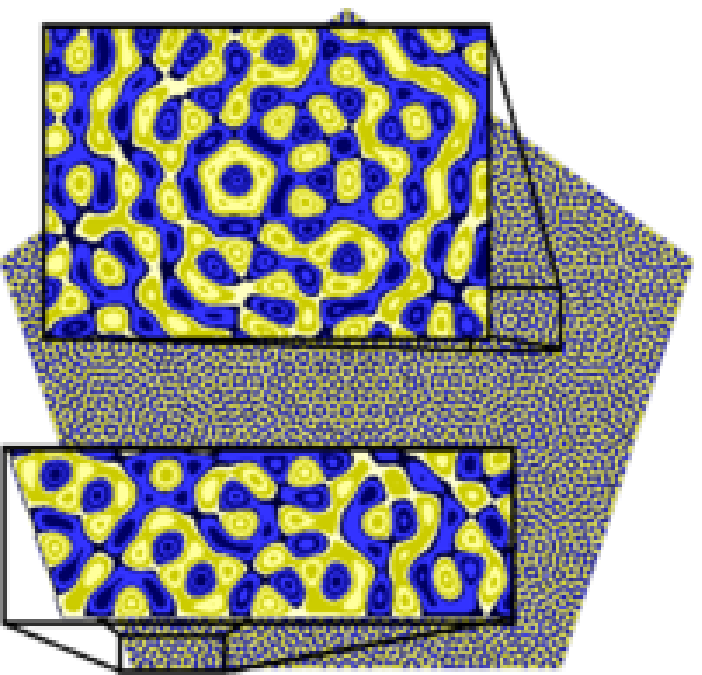}\label{fig:ds1396hb}}
\caption{Several representative Dirichlet eigenfunctions within the 
  unit-edged regular pentagon.  Blue and yellow areas are opposite
  sign, and light-dark steps indicate level contours. \mbox{(a)~The}
  degenerate pair just below \mbox{DS-18} \mbox{(see
    Fig.~\ref{fig:starandpentagon}, LEFT)}. \mbox{(b)~Another}
  example that strongly reveals the {\em pentagram\/} nodal
  pattern. \mbox{(c) The} first eigenfunction with eigenvalue
  above 100,000, corresponding to the 8139\textsuperscript{th} overall
  eigenvalue, and having about fifty wavelengths along a pentagon edge
  \mbox{($\Lambda\approx 0.020$)}.  Notice how the boundary shape is
  revealed within the sea of apparent chaos (upper inset).  See 
   Eq.~(\ref{eq:8139}) for an abbreviated 300-digit eigenvalue bound.}
\label{fig:pentagonplots}
\end{figure}

\begin{figure}[htb!]
\includegraphics[scale=0.7]{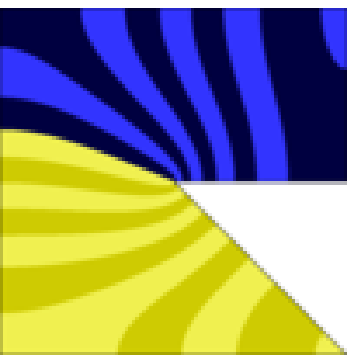}~~\includegraphics[scale=0.7]{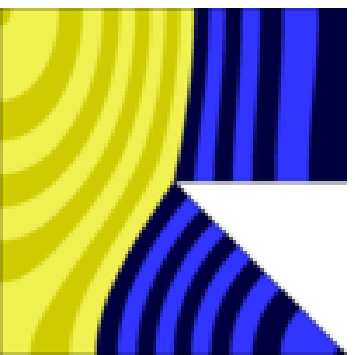}~~\includegraphics[scale=0.7]{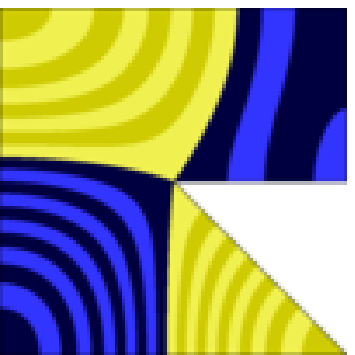}~~~~~\includegraphics[scale=0.7]{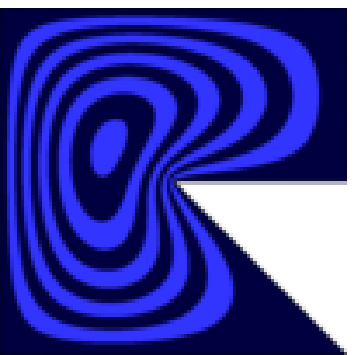}~~\includegraphics[scale=0.7]{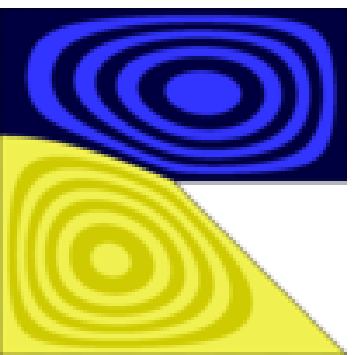}~~\includegraphics[scale=0.7]{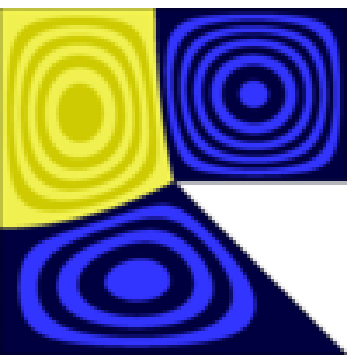}
\caption{Example contour plots.  Blue and yellow areas are opposite
  sign, and light-dark steps indicate level contours. These are (LEFT) the
  lowest three Neumann and (RIGHT) the lowest three Dirichlet
  eigenfunctions. Their eigenvalues are listed in
  Table~\ref{tab:ultraprecise}, and they are not degenerate.}
\label{fig:cutsquareeigenfunctions}
\end{figure}

It is generally a good idea to plot some eigenfunctions if only to
inspect the contours and nodal patterns to ensure that one is actually
calculating eigenvalues.  Several eigenfunction contour plots are
shown in Figs.~\ref{fig:stareigenfunctions}, \ref{fig:pentagonplots},
and~\ref{fig:cutsquareeigenfunctions}.

Of course, published concerns regarding the numerically
ill-conditioned nature of this method must be addressed. Some such
concerns were actually identified by FHM, but more recently clarified
by Betcke and Trefethen~\cite{bt2004} in 2004, when they demonstrated
the numerical advantages of the so-called ``generalized singular-value
decomposition'' (GSVD) method, a relative to the point-matching
method.

Fortunately, the answer to make the point-matching method work is
short and simple: One must both (a)~select adequate matching points
(both number and distribution) and (b)~keep enough precision in the
intermediate calculations.  My empirical observations are that
Chebyshev nodes chosen as matching points (as suggested by Betcke and
Trefethen~\cite{bt2004}) often work well, even where equally-spaced
points do not, and the precision of the intermediate calculations must
be significantly higher, often several times, than the precision of
the eigenvalue.

The ability to calculate many digits depends on the convergence rate,
and a good way to improve it is to exploit symmetries. Fortunately,
this exploitation has the important side-effects of (a)~reducing the
complexity of the region in which one must work, (b)~classifying
eigenmodes, and (c)~disentangling geometrically degenerate
eigenfunctions. Do not underestimate the importance of exploiting
symmetry.

At virtually every step in this project, only free
software\footnote{Typically per the
GPL \url{http://www.gnu.org/philosophy/free-sw.html}.}  running on
modern commodity computer hardware has been used, and with great
success.  The software of choice is the {\tt GP/PARI}
calculator~\cite{PARI2} using the GNU Multiple Precision Arithmetic
Library {\tt gmp}~\cite{gmp6} running on a laptop computer with a {\tt
GNU/linux} operating system and its countless ancillary programs.
(Some use was also made of {\tt maxima}~\cite{maxima} for the
occasional symbolic calculation.)  This programming environment
permits efficient numerical computations while ``natively'' retaining
up to several thousand digits of precision.

Several personal computers were used, but the best was a modern laptop
computer with \mbox{16\,GB RAM} and a \mbox{quad-core i7} processor
with eight threads.\footnote{At times, using up to \mbox{50\,GB} swap
space.} CPU times are reported using that laptop, and are unavoidably
approximate due to multitasking.

Before beginning, it should be made clear that this is a classical and
very well-known problem, worked on by many people over the last
two-hundred years---with many applications and results.  As such, I
shall limit the discussion to only those facts that are required to
reproduce and possibly extend the present calculations.  A
thirty-year-old, but still popular and relevant survey of the problem
was given by Kuttler and Sigillito~\cite{ks1984}. More practically,
active investigators, Barnett and Betcke have created {\tt
  MPSpack}~\cite{bb2010,bh2014} that helps bring sophisticated
eigenvalue calculations closer to the rest of us.

Despite the long history, except for the recent work of P.~Amore, et~al.~\cite{abfr2015}
, who incidentally make this same observation, I
am unaware of any published, non-closed-form eigenvalues accurate to
just beyond a dozen or so digits for {\em any\/} shape not related to
the closed-form solutions within the equilateral triangle, rectangle,
or circle (ellipse). All eigenvalue results in this report are 
likely unprecedented.

\section*{\label{sec:eigen}The eigenvalue problem}
Let $\mathbf{r}$ be a point in the plane described by either Cartesian
coordinates \mbox{``$(x,y)$''} or polar coordinates \mbox{``$(r,\theta)$''},
where \mbox{$x=r\cos\theta$} and \mbox{$y=r\sin\theta$}, and where notation
ambiguity is removed by context.
The two-dimensional Helmholtz equation is
\begin{equation}
    \Delta\Psi(k;\mathbf{r})+\lambda\Psi(k;\mathbf{r}) = 0
\label{eq:helmholtz1}
\end{equation}
where~$\Delta$ is the Laplacian and $k=\sqrt{\lambda}$ is the usual
``wavenumber''. Without a boundary, $\lambda$ (or $k$) is treated as a
continuous eigen-parameter, and $\Psi(k;\mathbf{r})$ may describe a
free wave with wavelength $\Lambda=2\pi/k$.  The ``interior''
Helmholtz eigenvalue problem is obtained by restricting $\mathbf{r}$
to the interior of a region (one of Fig.~\ref{fig:allshapes}), and
imposing relevant boundary conditions (Neumann or Dirichlet), which
constrains $\lambda$ to a non-accumulating set of discrete
eigenvalues, some of which may be degenerate. (The trivial,
``closed-form'', Neumann solution, i.e., $\Psi=\mbox{constant}$ with
$\lambda=0$, is completely ignored in this project.)

The point-matching method works if the eigenvalues are
non-degenerate. Thus it is very important to deal with degeneracies
either by (a)~dismissing them or (b)~disentangling them using
symmetry.

First, closed-form solutions are not only known\footnote{Among the
  chosen shapes, the L-shape, cut-square, and regular hexagon have
  subsets of closed-form solutions. For polygons, a closed-form
  solution can be written as a finite sum of plane waves, and its
  eigenvalue is related to~$\pi$, presently known to $13.3\times
  10^{12}$ digits.}, but also arbitrarily high in ``accidental''
degeneracy as one climbs the eigenvalue towers\footnote{This is
  related to the increasing number of integer solutions to Diophantine
  equations $a^2=b^2+c^2$ (rectangle) and $a^2=b^2+bc+c^2$
  (equilateral triangle).}; so after identifying them, simply exclude
them from the calculations.

Second, if there is a geometric or reflection symmetry, all
eigenfunctions can be sorted into separate symmetry classes---some of
which may be closed-form. Doing so effectively splits the
problem up into a set of sub-problems, each one with a corresponding
symmetry-reduced polygon and edge conditions, and a resulting,
non-accumulating, infinite tower of distinct eigenvalues
\begin{equation}
 0 < \lambda_1 < \lambda_2 < \lambda_3 < \cdots < \lambda_\alpha < \cdots,
\end{equation}
where $\alpha$ labels the eigenvalue within that tower. Considering
only the non-closed-form symmetry classes, it is indeed assumed
that there are no degeneracies within such a
tower.\footnote{Such non-closed-form degeneracies are
  possible but presumably quite rare. This assumption should be
  independently verified.}  Thus, to each non-closed-form
symmetry class, we associate a sub-problem, effectively consisting of
non-degenerate eigenmodes. A tower of such 
eigen-pairs for a given symmetry class shall be written as the set
\begin{equation}
       \left\{ \lambda_\alpha, \Psi_\alpha(\mathbf{r})\,\middle|\, \alpha=1,2,3,... \right\}
       \label{eq:exacteigenpair}
\end{equation}
where
$\Psi(k_\alpha;\mathbf{r})\equiv\Psi_\alpha(\mathbf{r})$. Degeneracies,
of course, may exist between separate symmetry classes, so it may be possible
to choose a subset of non-closed-form classes and calculate only those eigenvalues.

It is these symmetry-reduced, sub-problems that we focus on. In the
end, one can piece together all of the eigenmodes of each symmetry
class to construct the complete spectrum for the chosen shapes in
Fig.~\ref{fig:allshapes}, which restores the rich assortment of
symmetries and degeneracies. See, for example,
Fig.~\ref{fig:starSpectrum} depicting the lowest seventy-nine
eigenvalues of the star.

As is conventional, an edge or a line of symmetry within the shape is
said to be ``even'' or ``odd'' according to how a corresponding
eigenfunction transforms under Riemann-Schwarz
reflection.
Also, for the purposes of this report, the
analyticity of a vertex is defined by whether or not an eigenfunction
can be (locally) Riemann-Schwarz reflected around that vertex and
remain single valued.\footnote{All non-closed-form solutions in this
  project have an $\Omega$ with exactly one non-analytic vertex.  For
  all eigenfunctions within a given symmetry class, it seems that a
  vertex will be either analytic or not.}

In the current project, not all boundary edges will be either even or
odd. Indeed, types of symmetry arise where (a)~the eigenfunction
satisfies a periodic-like condition such that values along one edge
are proportional to values along another edge (star and regular
polygon), or (b)~a reflection symmetry may be exploited, effectively
excluding terms from the eigenfunction expansion (L-shape and
cut-square).

\begin{figure}[htb!]
\includegraphics[scale=0.8]{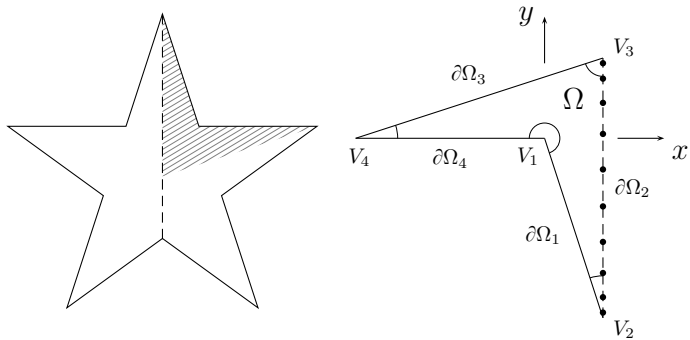}
\caption{(LEFT) Star shape showing how (RIGHT) an arrowhead polygon,
  $\Omega$, in its canonical position, is constructed for the
  $\mathcal{B}_e$ and $\mathcal{C}_e$ eigenfunctions of the star. The
  vertical bisector of the star for these eigenfunctions is an even
  (dashed) line, a portion of which becomes $\partial\Omega_2$. Edges
  $\partial\Omega_1$ and $\partial\Omega_4$ have same boundary
  conditions as the star, while edge $\partial\Omega_3$ satisfies a
  periodic-type boundary condition (values on $\partial\Omega_3$ are
  proportional to corresponding values on $\partial\Omega_2$). Also
  shown are $N/2=10$ Chebyshev distributed matching points on
  $\partial\Omega_2$. The other $N/2$ points would be symmetrically
  located on $\Omega_3$.  $V_1$ is non-analytic for all symmetry
  classes of the star.}
\label{fig:starA}
\end{figure}

To demonstrate the technique, it is necessary to introduce a polygon,
$\Omega$, as illustrated in the example of Fig.~\ref{fig:starA}. The
specific form of $\Omega$ depends on the symmetry group of the
problem, and for this project, define it to be an $s$-sided polygon
with at most one non-analytic vertex from which all edges can be seen.
Label the vertices $V_a$, where $a=1,2,...,s$; beginning with the
non-analytic vertex and proceeding clockwise. (In the example, $s=4$.)
The edges are labeled similarly, specifically, $\partial\Omega_a$
connects $V_a$ to $V_{a+1}$ (\mbox{mod~$s$}, as required). The
``adjacent'' edges are $\partial\Omega_1$ and $\partial\Omega_s$,
while the ``point-matched'' edges are $\partial\Omega_a$, where
$a=2,3,...,s-1$. In this project, the shapes are such that only the
adjacent edges lie along lines passing through the non-analytic
vertex.

Define the {\em canonical position\/} of $\Omega$ as follows.  First
locate the non-analytic vertex, $V_1$, at the origin, and orient as
needed so that $\partial\Omega_1$ and $\partial\Omega_2$ are
independently even or odd. Then, rotate around $V_1$ until
$\partial\Omega_2$ is parallel to the \mbox{$y$-axis} at a positive
\mbox{$x$-value}.  In this canonical position, each edge
$\partial\Omega_a$ is a segment of the line
\begin{equation}
 \mathcal{L}_a : \, A_a x+ B_a y = C_a 
\label{eq:linea1}
\end{equation}
where the constants ($A_a$, $B_a$, $C_a$) are uniquely chosen by
requiring that
\begin{equation}
  \widehat{\mathbf{r}}_{n,a} = A_a \widehat{\mathbf{x}} + B_a \widehat{\mathbf{y}}
\label{eq:linea2}
\end{equation}
be an outward pointing unit normal. 

The two adjacent edges, $\partial\Omega_1$ and $\partial\Omega_s$, at
polar angles $\phi_1$ and $\phi_s$, respectively, form the
non-analytic vertex, which has internal angle
\begin{equation}
            \Delta\phi=\phi_s-\phi_1.
\end{equation}
It is also convenient to introduce the abbreviation,
\begin{equation}
      \tilde\theta \equiv \theta - \phi_1
 \label{eq:thetaabbreviation}
\end{equation}

Having stated the problem, with some specific restrictions and
conventions, the next task is to describe the method of solution, the
so-called ``point-matching'' method.

\section*{The point-matching method}
Precisely following FHM, expand the eigenfunction in a truncated
Fourier-Bessel series
\begin{equation}
  \Psi^{[N]}(k;\mathbf{r}) = \sum_{\nu=1}^{N} c_\nu^{[N]}\, 
                \psi_\nu(k;\mathbf{r}) \vspace{1ex}
\label{eq:fourierbessel}
\end{equation}
where the so-called ``basis functions'' are 
\begin{equation}
   \psi_\nu(k;\mathbf{r}) = J_{m_\nu}(kr)\,\times\,
            \left\{
              \begin{array}{ll}
                    \sin(m_\nu\tilde{\theta}), & \mbox{\ $\partial\Omega_1$(odd)}  \vspace{1ex} \\
                    \cos(m_\nu\tilde{\theta}), & \mbox{\ $\partial\Omega_1$(even)}
              \end{array}
            \right.
\label{eq:basisfunctions}
\end{equation}
where, as indicated, the choice is based solely on if the adjacent
edge $\partial\Omega_1$ is odd or even; and $J_m(x)$ is the Bessel
function of the first kind of order~$m$.\footnote{To calculate only
  eigenvalues, we need not normalize these basis functions.}  The
index~$\nu$ runs from~$1$ to~$N$, here labeling the term in the
expansion; and the $m_\mu$ are not yet specified.

When the symmetry requirements and boundary conditions are enforced
along {\em both\/} adjacent edges, the \mbox{$m$-values} are
restricted.  It often happens that the adjacent edges are
independently even or odd, in which case, the $m$-values are given by
\begin{equation}
  m_\nu = \frac{\pi}{\Delta\phi} \times \left\{ \begin{array}{ll}
                      (\nu-1) & \mbox{ both even} \\
                      (\nu-\frac{1}{2}) & \mbox{ one odd, other even} \\
                      \hphantom{(}\nu & \mbox{ both odd} \\
                   \end{array}
           \right. 
\label{eq:mvalues1}
\end{equation}
Among the chosen shapes, the L-shape and cut-square have other
interesting symmetry-related restrictions on the $m$-values due to the
fact that these shapes can be obtained by reflecting a 45-90-45
triangle, respectively, six and seven times.

Once the required $m$-values are selected, it is important to realize
that $\Psi^{[N]}(k;\mathbf{r})$ is an {\em exact\/} solution within
the infinite sector \mbox{$\phi_1\le\theta\le\phi_2$}. This is true
for any positive integer $N$ and real $\lambda=k^2>0$ (excluding
$\lambda=0$). The function \mbox{$\Psi^{[N]}(k;\mathbf{r})$} is
specifically crafted to handle the non-analyticity of~\mbox{$V_1$} and
the boundary or symmetry conditions on the adjacent
edges~\mbox{$\partial\Omega_1$} and~\mbox{$\partial\Omega_s$},
exactly.

The above formalism generally applies to the MPS.  The next step,
unique to the ``point-matching'' method, is to select~$N$ points,
\begin{equation}
         \mathcal{P}_N: \left\{\mathbf{r_\mu} \,\middle|\,\mu=1,2,3,\cdots,N\right\}
\end{equation}
on the non-adjacent edges where the edge conditions may be
nontrivially enforced on those edges.\footnote{For example, if
  $\partial\Omega_1$ and $\partial\Omega_2$ are both odd, then $V_2$
  cannot be used as a matching point.} This is the same~$N$ that
appears in the eigenfunction expansion, Eq.~(\ref{eq:fourierbessel}).
Selecting adequate sets of matching points is a nontrivial task that
is key to making the point-matching method work.  For the immediate
following, assume this has been done.

Enforcing the boundary conditions or symmetry-related edge conditions
at the~$N$ matching points yields~$N$ linear equations in the~$N$
expansion coefficients (the $c_\nu^{[N]}$ in Eq.~(\ref{eq:fourierbessel})), and
a resulting \mbox{$N\times N$} ``point-matching matrix'',
$\mathcal{M}^{[N]}(\lambda)$.
Linear algebra suggests that the determinant
of this matrix must be zero for nontrivial solutions to exist, i.e.,
\begin{equation}
          \left|\mathcal{M}^{[N]}(\lambda)\right| = 0
\label{eq:pointmatchingdeterminant}
\end{equation}
This is the so-called ``point-matching determinant'', and finding
roots of this equation for increasing values of $N$ is the name of the
game.

For a given $N$, there is an infinite
number of roots of the point-matching determinant, but only the first
$\kappa_N$ of them,
\begin{equation}
     0 < \lambda_1^{[N]} < \lambda_2^{[N]} < \cdots < \lambda_{\kappa_N}^{[N]}
\end{equation}
are guaranteed to correspond to actual eigenvalues. 
The number $\kappa_N$ is such that if $N$ is increased to \mbox{$N+\Delta
  N$}, this set of roots may become better estimates of the
eigenvalues, but new roots may be introduced only for values of
$\lambda$ a little bigger than~\mbox{$\lambda_{\kappa_N}^{[N]}$}.  

Thus, for a given eigenvalue, $\lambda_\alpha$, there is a minimum
$N$-value, say $N_{\alpha,1}$, above which
Eq.~(\ref{eq:pointmatchingdeterminant}) may yield increasingly better
estimates as $N$ is incremented. I conjecture that a more detailed
picture shall require $N_{\alpha,1}$ be big enough so that a gently
``wiggling'' contour passes continuously through all the matching
points---a contour along which the edge conditions are satisfied
exactly.\footnote{For example, if the point-matched edge is odd, then
  that contour would be a nodal curve passing through the matching
  points.  If $N$ is too small for a given eigenvalue, then, for
  example, two adjacent matching points might not be directly joined
  by nodal curve, although nodal curves will pass through those
  matching points.}

The smallest, problem-dependent value of $\Delta N$ defines a
``properly incremented'' set of $N$-values, say $\{N_{\alpha}\}$, where
its members are
\begin{equation}
     N_{\alpha,i} = N_{\alpha,1}+(i-1)\Delta N \qquad i=1,2,3,\cdots;
\label{eq:properincrement}
\end{equation}
that yield a corresponding set of approximate eigenvalues
\begin{equation}
    \left\{ \lambda^{[N_{\alpha,1}]}_\alpha,\lambda^{[N_{\alpha,2}]}_\alpha,\lambda^{[N_{\alpha,3}]}_\alpha,\cdots\right\}.
\label{eq:approximateeigenvaluesequence}
\end{equation}

Provided everything ``works'', we then tacitly assume that the
approximate eigen-pair converges to the corresponding exact
eigen-pair (see Eq.~(\ref{eq:exacteigenpair})), specifically,
\begin{equation}
  \lim_{N\to\infty} \left\{\lambda^{[N]}_{\alpha},\Psi^{[N]}(k_\alpha^{[N]};\mathbf{r})\right\}
    = \left\{ \lambda_\alpha, \Psi_\alpha(\mathbf{r}) \right\}
\label{eq:eigenpairlimit}
\end{equation}
where $N\in\{N_\alpha\}$. 
The procedure depends on several factors, most
notably, the nature of the vertices and the distribution of the
matching points. 

\section*{Eigenvalue bounds}
If the rate of convergence is ``exponential'', a very interesting
thing happens that might be used to obtain eigenvalue bounds, by
inspection.  Specifically, the problem may be set up so that that
sequence of approximate eigenvalues,
Eq.~(\ref{eq:approximateeigenvaluesequence}), as $N$ increases, can be
made to alternate above and below an asymptote---assumed to correspond
to an eigenvalue---in such a way that the peaks of that alternation
effectively provide an increasingly narrow bound for that
eigenvalue.\footnote{I first observed this alternation property in
1993~\cite{phdthesis} with exponentially-convergent closed-form
solutions.}  This means we can write
\begin{equation}
      \lambda_\alpha^{[N_\downarrow]} < \lambda_\alpha < \lambda_\alpha^{[N_\uparrow]}
\label{eq:eigenvaluebound}
\end{equation}
where $N_\downarrow$ and $N_\uparrow$ correspond to \mbox{$N$-values}
of a given minimum and the previous or next maximum, respectively.
Table~\ref{tab:Lshapeloworder} clearly illustrates that alternation
for the L-shape calculation---not only for my modern calculation, but
(interestingly) also for the original FHM data.

To calculate many digits, one must estimate the roots of the
point-matching determinant, Eq.~(\ref{eq:pointmatchingdeterminant}),
to a precision somewhat higher than the observed bound,
Eq.~(\ref{eq:eigenvaluebound}), at a given \mbox{$N$-value}. By
construction, the eigenvalues are not degenerate, meaning the
point-matching determinant only has simple roots, so it is actually
quite straightforward to calculate those roots to very high
precision.\footnote{The point-matching determinant changes sign as one
passes a root, i.e., provided that $\delta\lambda$ is sufficiently
small, then
if \mbox{``$|\mathcal{M}(\lambda)|\times|\mathcal{M}(\lambda+\delta\lambda)|<0$''},
then the root is in the interval $(\lambda,\lambda+\delta\lambda)$.
Without knowing the smallest gap between eigenvalues, it is important
to independently locate each eigenvalue. The GSVD method is perhaps
better suited to sweep an interval to locate eigenvalues.}  This key
ingredient---the ability to precisely locate those roots---is a
significant advantage of the point-matching method over some other
methods, like the GSVD method (which seeks minima of a function).

Another important feature of the point-matching method is that the
number of basis functions is maximized for a given number of matching
points. With more terms, the approximate eigenfunction better matches
the exact eigenfunction. 

With a given eigenvalue bound in hand, there are several related and
useful numbers one can calculate. They are \mbox{(a)~the} relative
gap, \mbox{(b)~the} approximate number of correct digits, and
\mbox{(c)~the} global convergence rate; respectively,
\begin{equation}
\mbox{(a)}\quad\epsilon_\alpha^{[\widehat{N}]} = \frac{\lambda^{[N_\uparrow]}_\alpha-\lambda^{[N_\downarrow]}_\alpha}
                       {\frac12\left[{\lambda^{[ N_\uparrow ]}_\alpha+\lambda^{[N_\downarrow]}_\alpha}\right]}
\qquad\qquad\mbox{(b)}\quad   D_\alpha^{[\widehat{N}]} = -\log_{10}\left(\epsilon_\alpha^{[\widehat{N}]}\right)
\qquad\qquad\mbox{(c)}\quad  \rho^{[\widehat{N}]}_\alpha = \frac{D_\alpha^{[\widehat{N}]}}{\widehat{N}}
\label{eq:relativegap}
\end{equation}
where the $\widehat{N}=\max(N_\uparrow,N_\downarrow)$ form a new set
of $N$-values, \mbox{$\{\widehat{N}_\alpha\}$}, a subset of $\{N_\alpha\}$.

Note that Eqs.~\mbox{(\ref{eq:relativegap})} provide a practical,
numerical definition of what it means to be ``exponentially
convergent''. Such a practical definition is useful since the exact
value of the eigenvalue is never really known. Specifically, if
$\rho_\alpha^{[\widehat{N}]}=\rho_\alpha$ is constant
wrt~$\widehat{N}$, then the eigenvalue {\em relative gap\/} decays
exponentially, i.e.,
\begin{equation}
 \epsilon_\alpha^{[\widehat{N}]} \sim 10^{- \rho_\alpha  \widehat{N}}
\end{equation}
as $\widehat{N}\to\infty$.  My observation is that the global convergence
rate~$\rho$ for non-closed-form solutions is typically a deceasing
function of~$\widehat{N}$.

Convergence rates that permit relatively easy \mbox{100-digit} results
range down to about~$1/6$ or so, and the best convergence rates for
the chosen shapes are typically a little better than~$1/2$.  Note that
this reciprocal notation is useful because it reveals how many
matching points---or terms in the expansion---must be added to achieve
each additional eigenvalue digit: Thus, if $\rho=1/4$, then four
additional terms in the expansion and a corresponding four additional
matching points will add an additional eigenvalue digit.

To report a bound, suffix the matching digits with a rounded-up
superscript and a rounded-down subscript. This conservative approach
helps ensure that the reported bound will include the true
eigenvalue. To illustrate, consider the relatively
difficult-to-calculate, lowest Dirichlet eigenvalue within a 256-sided
regular polygon (area=$\pi$). It is first precisely calculated
to
\begin{equation}
 \lambda=\left\{ \begin{array}{ll}
  5.7831876203689428757 2892661 ... & \mbox{ for $N_\downarrow= 930$} \\
  5.7831876203689428757 8671658 ... & \mbox{ for $N_\uparrow=   929$} \\
         \end{array}\right.
\end{equation}
This calculation used $N=282, 283,...,929,930$; while judiciously
skipping $N$-values, and stopping when
$\epsilon<10^{-20}$. Very reliably, even and odd $N$-values provided
lower and upper bounds, respectively. From that result, one may write
the bound
\begin{equation}
      \lambda=  5.7831876203689428757\,_{289}^{868}
\label{eq:polygon256}
\end{equation}
for $\widehat{N}=930$. For this example, \mbox{$\epsilon=9.99\times 10^{-21}$},
\mbox{$D=20.00$}, and~\mbox{$\rho=1/46.5$}. To get that twenty-digit
result took almost two CPU days.  Despite achieving {\em only\/} twenty digits,
this example does exhibit exponential convergence. Indeed, more
detailed inspection reveals the relationship
\begin{equation}
  \epsilon^{[\widehat{N}]}\approx 10^{-\left(5.73+ \widehat{N}/65.1\right)}
\end{equation}
based on $280\le\widehat{N}\le 930$, which specifically indicates
exponential convergence.\footnote{Going from $\widehat{N}=280\to
  930$, the global convergence rate decreases from $1/27.9\to
  1/46.5$. Asymptotically, $\rho^{[\widehat{N}]}\sim 1/65.1$.}
Incidentally, for comparison, the best ``published'' value may be
calculated with Eq.~(\ref{eq:regpolasymp}), yielding
\mbox{$\lambda\approx \underline{5.78318762036894}6869...$.}, of which the first
fifteen (underlined) digits appear to be correct. I am unaware of any
efforts to calculate this particular eigenvalue.

Note that above a few dozen eigenvalue digits, either the ``correct
digits'' \mbox{(e.g., $5.78...28757$)} or the ``correctly rounded digits''
\mbox{(e.g., $5.78...2876$)} may seem more impressive than the actual bound.
But the bound is nevertheless useful because the exact value is not known, and
one may not know otherwise where to stop counting correct digits.

It is the alternation property that provides the bounds, but what
causes that alternation?  By considering the necessary continuous
contour through the matching points along which the boundary or
symmetry conditions are satisfied exactly, I offer a simple heuristic
explanation.\footnote{This is certainly not a mathematical proof, but
  is quite plausible. } Such a contour defines an area $\mathrm{A}_N$
in which $\{\lambda^{[N]}_\alpha,\Psi^{[N]}(k_\alpha,\mathbf{r})\}$ is
an {\em exact\/} solution. As matching points are added, that area not
only approaches the polygon area, $\mathrm{A}$, but I conjecture that
the contour ``flips'', analogous to ``$\sin\to-\sin$''; and that such
a flipping causes the area to alternate above and below the polygon
area.  Since the area difference is very small\footnote{The area
  difference for a 100-digit eigenvalue is $\Delta A/A=10^{-100}$:
  Whimsically compare the area of a drum (vibrating membrane) with a
  diameter matching that of the universe, i.e., $R=47\times
  10^9 \,\ell \mathrm{y}$, $A=6.1\times 10^{53}\,\mathrm{m}^2$, to the
  area change equal to the cross-sectional area of a proton, i.e.,
  $R=0.88\,\mathrm{fm}$, $\Delta A = 2.4\times
  10^{-30}\,\mathrm{m}^2$. This $\Delta A/A=4.0\times 10^{-84}$ is some
  fifteen orders of magnitude larger.}  , using
$\lambda^{[N]}_\alpha\mathrm{A}_N\approx \lambda_\alpha\mathrm{A}$, it
becomes plausible that such a flipping causes the alternation.

My empirical observation is that if there is exponential convergence,
then it becomes more likely that flipping occurs with each proper
increment the closer the matching-point distribution is to being
equally-spaced. More optimal (or necessary) Chebyshev
distributions give better convergence rates, but quite often require
several proper increments to flip. In that case, the number of proper
increments is often irregular, but may become more regular for higher
\mbox{$N$-values}. Any automatic program should make sure that
what appear to be upper and lower limits are actually so.

Note that this eigenvalue bounding technique is related to, but
slightly different from the conventional approach of using the
approximate eigenfunction along the point-matched edges, as in the
Moler-Payne method~\cite{mp1969}.  A practical difference is that one
need not calculate the coefficients (to evaluate the eigenfunction
along the point-matched edges) to calculate the upper and lower limits
to the eigenvalue. Instead, the eigenvalue bound here is identified by
``simply'' watching roots of the point-matching determinant as $N$ is
increased. However, ease of calculation may not be the only advantage:
It may indeed provide {\em better\/} bounds, as demonstrated with
the L-shape calculation below.

To make all of this work well, one must identify symmetries, select
polygons and adequate matching point distributions, and calculate
precise values of the point-matching determinant for a properly
incremented set of $N$-values. For a given problem, this may be an
iterative procedure as one discovers new things.

\section*{The matching points}
Next consider the nontrivial task of choosing matching points.  Each
problem will have a specific choice, but there are some common
strategies.  Suppose each point-matched edge $\partial\Omega_a$ has
$n_a$ matching points at which the edge conditions may be nontrivially
satisfied.

For the point-matching method to work well, it is somewhat important
that---as $N$ is incremented---the same proportion of matching points
be used on each edge, and the spacing and distribution should also be
approximately similar on each edge.  These requirements help define a
``proper increment'', Eq.~(\ref{eq:properincrement}). 

It is also recommended that the maximum gap between matching points is
always smaller than some empirically determined fraction of the free
wavelength, \mbox{$\Lambda=2\pi/k$}. This will help ensure that the
contour connecting matching points does indeed pass
continuously through all the matching points.

Equally-spaced matching-points seem quite popular and easy to program,
but these do not always work and are certainly not the best.  It seems
that crowding points near vertices is good practice since that is
usually where the contour tends to deviate from the polygon edge the
most: Crowding the points constrains the contour better.  

Chebyshev nodes usually yield excellent results and happen to be quite
easy to program.  Indeed, let~\mbox{$s=0$} to~$\ell$ measure the
position along a point-matched edge $\partial\Omega_a$, then
\begin{equation}
   s_\mu = \frac{\ell}{2}\left[1 - \cos\left(\frac{\mu-\frac{1}{2}}{n_a}\,\pi\right)\right]
\label{eq:canonicalchebyshevnodes}
\end{equation}
where $\mu=1,2,3,...,n_a$ is a canonical set of Chebyshev nodes. This
distribution has points very close to the end-points, and a relative
crowding near those end-points. Variations may be used, for example,
to ensure more crowding near acute, distant vertices.

\section*{The point-matching matrix}
Armed with an adequate selection of matching points, the next task is
to describe the matrix elements. Each point-matched edge,
$\partial\Omega_a$ with $n_a$ matching points, can be odd, even, or
satisfy a periodic-type requirement; and each will yield $n_a$ rows of
the point-matching matrix.

\subsection*{Odd point-matched edge}

If the point-matched edge $\partial\Omega_a$ is odd, the $n_a$ rows
are obtained by equating \mbox{$\Psi^{[N]}(k;\mathbf{r})$} at the
$n_a$ matching points to zero. Specifically,
\begin{equation}
         \Psi^{[N]}(k;\mathbf{r_\mu}) = 0 \qquad \mbox{$\mathbf{r_\mu}\in\partial\Omega_a$(odd)},
\end{equation}
which yields $n_a$ rows,
\begin{equation}
         \mathcal{M}^{[N]}_{\mu\nu}(\lambda) = J_{m_\nu}(kr_\mu)\,\times\,
            \left\{
              \begin{array}{l}
                    \sin(m_\nu\tilde{\theta}_\mu),\mbox{ $\partial\Omega_1$(odd)}  \vspace{1ex} \\
                    \cos(m_\nu\tilde{\theta}_\mu),\mbox{ $\partial\Omega_1$(even)}
              \end{array}
            \right.
\label{eq:matrixelementsodd}
\end{equation}
The vast majority of published accounts of all variants of the MPS
use Dirichlet boundary conditions, so this is perhaps a well-known
result. It is also relatively easy to program.

\subsection*{Even point-matched edge}
If the point-matched edge $\partial\Omega_a$ is even, the formulas are
also elementary and easy to program, but a little long. These results
appear to be less well-known, so they may be of some utility.

First recall equations~(\ref{eq:linea1}) and~(\ref{eq:linea2}), which
define a \mbox{line---a} segment of which forms
\mbox{edge~$\partial\Omega_a$---and} an ``outward-pointing'' unit normal,
$\widehat{\mathbf{r}}_{n,a}$.  Equating the normal derivative of
\mbox{$\Psi^{[N]}(k;\mathbf{r})$} at the $n_a$ matching points to zero
shall yield the corresponding $n_a$ rows of the matrix.  Specifically,
\begin{equation}
    \left. \frac{\partial\Psi^{[N]}(k;\mathbf{r})}{\partial r_{n,a}}\right|_{\mathbf{r}_\mu} = 0
   \qquad \mbox{$\mathbf{r_\mu}\in\partial\Omega_a$(even)}
\end{equation}
where
\noindent
\begin{equation}
\displaystyle
        \frac{\partial\Psi}{\partial r_{n,a}} = A_a\,\frac{\partial\Psi}{\partial x} + B_a\, \frac{\partial\Psi}{\partial y}.
\end{equation}

To make the formulas less cumbersome, split $\mathcal{M}$ into two
$N\times N$ matrices, $\mathcal{X}$ and $\mathcal{Y}$, such that
\begin{equation}
         \mathcal{M}^{[N]}_{\mu\nu}=A_a\mathcal{X}^{[N]}_{\mu\nu}+B_a \mathcal{Y}^{[N]}_{\mu\nu}.
\end{equation}
Differentiating the basis functions and evaluating at the matching points yields
\begin{widetext}
\begin{subequations}
\begin{eqnarray}
\displaystyle
 r_\mu^2\mathcal{X}^{[N]}_{\mu\nu}(\lambda) &=& 
         \left\{ \begin{array}{ll}\displaystyle
   m_\nu \left[
                       x_\mu \sin(m_\nu\tilde{\theta}_\mu) - y_\mu \cos(m_\nu\tilde{\theta}_\mu)
        \right] \, J_{m_\nu}(kr_\mu)
             -kr_\mu x_\mu \sin(m_\nu\tilde{\theta}_\mu) \, J_{1+m_\nu}(kr_\mu)
                        , & \mbox{$\partial\Omega_1$(odd)}\vspace{1ex} \\
\displaystyle
   m_\nu \left[
                       x_\mu \cos(m_\nu\tilde{\theta}_\mu) + y_\mu \sin(m_\nu\tilde{\theta}_\mu)
        \right] \, J_{m_\nu}(kr_\mu)
                                        -kr_\mu x_\mu\cos(m_\nu\tilde{\theta}_\mu) \, J_{1+m_\nu}(kr_\mu)
                ,  &\mbox{$\partial\Omega_1$(even)}
         \end{array}\right. \\
 r_\mu^2\mathcal{Y}^{[N]}_{\mu\nu}(\lambda) &=& 
         \left\{ \begin{array}{ll}\displaystyle
   m_\nu \left[
                       x_\mu \cos(m_\nu\tilde{\theta}_\mu) + y_\mu \sin(m_\nu\tilde{\theta}_\mu)
        \right] \, J_{m_\nu}(kr_\mu)
             -kr_\mu y_\mu \sin(m_\nu\tilde{\theta}_\mu) \, J_{1+m_\nu}(kr_\mu),
                 & \mbox{$\partial\Omega_1$(odd)}\vspace{1ex} \\
\displaystyle
   m_\nu \left[
                       -x_\mu \sin(m_\nu\tilde{\theta}_\mu) + y_\mu \cos(m_\nu\tilde{\theta}_\mu)
        \right] \, J_{m_\nu}(kr_\mu)
             -kr_\mu y_\mu\cos(m_\nu\tilde{\theta}_\mu) \, J_{1+m_\nu}(kr_\mu),&\mbox{$\partial\Omega_1$(even)}
         \end{array}\right.\qquad
\label{eq:matrixelementseven}
\end{eqnarray}
\end{subequations}
\end{widetext}
To calculate eigenvalues, one may absorb the factor $r_\mu^{-2}$ into
the expansion coefficients since it is common to each term in the row.

In most of my examples, $A_2=1$ and $B_2=0$, but Neumann eigenvalues
of the cut square shall also require $A_3=0$ and
$B_3=1$.

\subsection*{Periodic-type edge (dihedral symmetry)}
Next consider the periodic-type boundary conditions that arise if the
geometry has dihedral symmetry. Although this analysis is rather
elementary, it too seems to be rarely discussed in the context of this
problem. Since the regular polygon\footnote{Although most of this
  analysis applies to all regular polygons, I shall exclude the
  closed-form solutions, as explained above.} has dihedral symmetry,
it will be used to most efficiently develop the matrix elements; but
it also applies to another chosen shape, the 5-point star.

\begin{figure}[tb]
\includegraphics[scale=0.8]{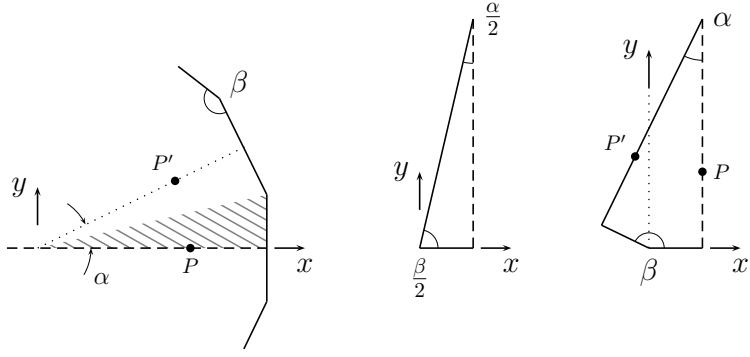}
\caption{(LEFT) A regular polygon with an even apothem (dashed) on the
  positive $x$-axis, (CENTER) the triangular fundamental region, and
  (RIGHT) the kite-shaped region for the doubly degenerate eigenmodes.
  Values of an even \mbox{($Q$-parity)}, doubly-degenerate
  eigenfunction are proportional to each other at points~$P$
  and~$P'$.}
\label{fig:regularpolygons}
\end{figure}

To that end, consider a $\sigma$-sided regular polygon (and other
diagrams) shown in Fig.~\ref{fig:regularpolygons}, with vertex angle
$\beta=(\sigma-2)\pi/\sigma$. The symmetry group of this shape is the
dihedral group, $D_\sigma$, of degree $\sigma$ or order $2\sigma$.  To
work out the symmetry properties, first center the polygon at the
origin with an apothem on the positive \mbox{$x$-axis}, as shown in
\mbox{Fig.~\ref{fig:regularpolygons}\,(LEFT)}.  The dihedral group is
generated by $R$, a counter-clockwise rotation by \mbox{$\alpha=
  2\pi/\sigma$}, and $Q$, a reflection through the \mbox{$x$-axis},
i.e., \mbox{$Q:y\to -y$}.

The dihedral group has only \mbox{1-dim} and \mbox{2-dim} irreducible
representations (irreps), so let $\eta_1$ and $\eta_2$ count those
irreps. If $\sigma$ is odd or even, then $\eta_1=2$ or $\eta_1=4$,
respectively.  Knowing $\eta_1$, we have
\mbox{$\eta_2=(2\sigma-\eta_1)/4$}.  These 1-dim and 2-dim irreps lead
to, respectively, non-degenerate and doubly-degenerate towers of
eigenvalues.

Since $Q^2=1$, define a \mbox{``$Q$-parity''} such that any
non-degenerate eigenfunction is either even or odd according to that
parity. For the doubly degenerate pairs of eigenfunctions, one can
always be made odd and the other even, and assume that this is done.
Use a subscript to denote the \mbox{$Q$-parity}, as in
\begin{subequations}
\begin{eqnarray}
   \varphi_e(x,-y)=+\varphi_e(x,y) \qquad \varphi_o(x,-y)=-\varphi_o(x,y)\qquad &&\mbox{Cartesian}\vspace{1ex}\\
   \varphi_e(r,-\theta)=+\varphi_e(r,\theta) \qquad \varphi_o(r,-\theta)=-\varphi_o(r,\theta)\qquad &&\mbox{Polar}
\end{eqnarray}
\end{subequations}
where ``$\varphi(\mathbf{r})$'' denotes a function in the regular
polygon centered at the origin.  An obvious but relevant fact is that
all odd $Q$-parity functions are zero on the \mbox{$x$-axis}, i.e.,
\mbox{$\varphi_o(x,0)\equiv 0$}.

For the regular polygon, the fundamental region may be chosen to be
the shaded triangle of Fig.~\ref{fig:regularpolygons}\,(LEFT) because
the entire polygon may be obtained from that triangle via group
operations. Relevant to this project, it is a right triangle with two
other angles
\begin{equation}
  \frac{\alpha}{2}=\frac{\pi}{\sigma}\qquad\mbox{and}\qquad    \frac{\beta}{2} = \frac{(\sigma-2)\pi}{2\,\sigma},
\end{equation}
exactly one of which, $\beta/2$, is non-analytic (except for
$\sigma=3$, $4$, and sometimes~$6$).

For the non-degenerate eigenfunctions, the triangular fundamental
region becomes~$\Omega$.  If $\sigma$ is even, there are eight
possible sets of edge conditions on this triangle since all its edges
can be independently even or odd.  But, if $\sigma$ is odd, there are
only four possible sets because the regular polygon's apothem must
have the same symmetry as its circumradius (line connecting its center
to a vertex).\footnote{This counting and the value of $\eta_1$ are
  related.}  When $\Omega$ is put in its canonical position,
\mbox{Fig.~\ref{fig:regularpolygons}\,(CENTER)}, it should become
clear how to calculate the matrix elements for the non-degenerate
eigenvalues.

For the doubly-degenerate eigenfunctions, we have an interesting
choice.  Since the group transformations form linear combinations of
degenerate eigenfunctions, we may either (a)~solve for {\em both\/} of
the degenerate eigenfunctions within the fundamental region or
(b)~solve for {\em one\/} of the degenerate eigenfunctions in an area
twice as large. In both cases, we can reconstruct both functions using
the $2\sigma$ group transformations. Since choice~(b) requires 
only one function, it shall be the better choice.

By reflecting the fundamental region about its hypotenuse to form the
kite-shaped quadrilateral, we obtain the only polygon that is both
twice as large and has (at most) one non-analytic vertex from which all
edges can be seen. This kite quadrilateral thus becomes~$\Omega$ for
the doubly-degenerate eigenmodes.\footnote{The other obvious choice,
the \mbox{$\alpha$-$\beta/2$-$\beta/2$} triangle, has {\em two\/}
non-analytic vertices.}

Cut this kite-quadrilateral out of the polygon as shown in
\mbox{Fig.~\ref{fig:regularpolygons}\,(LEFT)}, but before re-orienting
it, observe that if we ``rotate'' the {\em even\/} $Q$-parity
eigenfunction, we get the simple but important result that
\begin{equation}
   \varphi_e(r,\alpha)=\cos(\gamma\alpha)\,\varphi_e(r,0)
\label{eq:periodicphi}
\end{equation}
where $\gamma=1,2,...,\eta_2$. The effective purpose of $\gamma$ is to
identify to which one of the $\eta_2$ doubly-degenerate eigenvalue
towers this eigenfunction belongs.\footnote{To unify the expressions,
  one may include $\gamma=0$ (all~$\sigma$) and $\gamma=\sigma/2$
  (even~$\sigma$) for the even-parity, non-degenerate towers, but this
  is not done here because $P'$ is not on the boundary of triangular
  fundamental region. It is easier to keep the solutions corresponding
  to \mbox{1-dim} and \mbox{2-dim} irreps separated.} (Note that this
is where we used \mbox{$\varphi_o(r,0)\equiv 0$}).

That periodic-type symmetry relationship, Eq.~(\ref{eq:periodicphi}),
is useful because it relates values of an (even) eigenfunction on the
positive $x$-axis (an apothem) to values of the same (even)
eigenfunction on a neighboring apothem, i.e., at points~$P$ and~$P'$.

When that kite quadrilateral is re-oriented to its canonical position,
\mbox{Fig.~\ref{fig:regularpolygons}\,(RIGHT)}, the apothems become
point-matched edges and the regular-polygon boundary segments become
adjacent edges; so that Eq.~(\ref{eq:periodicphi}) gets transformed
into
\begin{equation}
   \Psi_e(k;r,\beta-\theta)=\cos(\gamma\alpha)\,\Psi_e(k;r,\theta),
\label{eq:periodicpsi}
\end{equation}
where the even apothem becomes $\partial\Omega_2$. Note how the points
$P$ and $P'$ are mapped.  The subscript ``$e$'' is retained to remind
us that we only need to consider this member of the degenerate pair to
determine the eigenvalue tower.

Since the regular polygon is assumed to have either Dirichlet or
Neumann boundary conditions, the adjacent edges
($\partial\Omega_1$~and~$\partial\Omega_4$) are either both odd or
both even.  In this case, referring to Eq.~(\ref{eq:mvalues1}), we can
use
\begin{subequations}
\begin{eqnarray}
       && \sin(\nu\,(\beta-\theta)) = -(-1)^\nu \sin(\nu\,\theta) \\
       && \cos((\nu-1)(\beta-\theta)) = -(-1)^\nu \cos((\nu-1)\,\theta)
\end{eqnarray}
\end{subequations}
to combine Eqs.~(\ref{eq:fourierbessel}) and (\ref{eq:periodicpsi}), yielding,
\begin{equation}
 0 = \sum_{\nu=1}^N c_\nu^{[N]}\, \left\{(-1)^\nu +  \cos(\gamma\alpha)\right\} 
      \, J_{m_\nu}(kr)
                    \times\left\{\begin{array}{ll}
                             \sin(m_\nu\theta),&\mbox{$\partial\Omega_1$(odd) and $\partial\Omega_4$(odd)} \vspace{0.5ex}\\
                             \cos(m_\nu\theta),&\mbox{$\partial\Omega_1$(even) and $\partial\Omega_4$(even)} 
                                        \end{array}\right.
\label{eq:evensymmetrycondition}
\end{equation}
Evaluating this equation at each of the matching points on $\partial\Omega_2$ establishes
$N/2$ equations for the $N$ expansion coefficients; while the other $N/2$ equations
are obtained by using the fact that $\partial\Omega_2$ is even. Requiring 
nontrivial solutions (for the expansion coefficients) leads to 
these $N/2$ rows of the point-matching matrix
\begin{equation}
  \mathcal{M}^{[N]}_{\mu\nu}(\lambda) = \left\{(-1)^\nu + \cos(\gamma\alpha) \right\} 
     \, J_{m_\nu}(kr_\mu)
                    \times\left\{\begin{array}{ll}
                             \sin(m_\nu\theta_\mu),&\mbox{$\partial\Omega_1$(odd) and $\partial\Omega_4$(odd)} \vspace{0.5ex}\\
                             \cos(m_\nu\theta_\mu),&\mbox{$\partial\Omega_1$(even) and $\partial\Omega_4$(even)} 
                                        \end{array}\right.
\label{eq:matrixelementsmixed}
\end{equation}
while the other $N/2$ rows are obtained using the fact that $\partial\Omega_2$ is even.

There are really two point-matched edges, but that periodic-type
boundary condition, Eq.~(\ref{eq:periodicpsi}), enables us to use one
point-matched edge, twice. This requires that $N$ be even, with $N/2$
points on edge $\partial\Omega_2$, and $N/2$ corresponding ``phantom''
points on $\partial\Omega_3$ which aren't actually numerically needed. 

As it happens, the symmetry group of the star is $D_5$, i.e., the same
as the regular pentagon.  The fundamental region is a triangle with
one non-analytic vertex (see Fig.~\ref{fig:starB}), and the region
$\Omega$ for the doubly degenerate eigenfunctions is an arrowhead
quadrilateral (see Fig.~\ref{fig:starA}), also with one non-analytic
vertex---this time, a re-entrant vertex. To accommodate the star (and
other polygons with dihedral symmetry),
Eq.~(\ref{eq:matrixelementsmixed}) is changed by replacing
$\theta_\mu$ with \mbox{$\tilde\theta_\mu=\theta_\mu-\phi_1$}, see Eq.~(\ref{eq:thetaabbreviation}).

As is common practice, the eigenfunctions can be identified by their
nodal and antinodal patterns. With dihedral symmetry, the
non-degenerate eigenfunctions have a maximal set of criss-crossing
even and odd lines of symmetry. These are obtained by Riemann-Schwarz reflecting
the eigenfunction within the triangular fundamental region to flesh out the
full eigenfunction. 

Every problem with such dihedral symmetry\footnote{Not just the
  regular polygons or the star. It is an exercise in geometry to
  figure out all of the specific shapes with dihedral-symmetry that
  also yield exponentially-convergent eigenfunction expansions (per
  the current procedure).}  will have at least two towers of
non-degenerate eigenfunctions, here named the ``symmetric''
$\mathcal{S}$ and the ``antisymmetric'' $\mathcal{A}$, which have all
even and odd lines of symmetry, respectively. If $\sigma$, as in
``$D_\sigma$'', is even, then there are two more non-degenerate
symmetry classes that have lines of symmetry that alternate, even and
odd, as one proceeds around the polygon. These might be named
$\mathcal{S}'$ and $\mathcal{A}'$ according to even and odd
$Q$-parity, respectively.

The doubly degenerate eigenfunctions can be similarly
identified by the nodal and antinodal patterns since these functions
have [at least] $\gamma$ nodal {\em curves\/} crossing at the
origin. Since the complexity of the nodal pattern increases with
$\gamma=1,2,...,\eta_2$; and we already used the symbol $\mathcal{A}$, label
these $\eta_2$ symmetry classes using $\mathcal{B}$, $\mathcal{C}$,
$\mathcal{D}$, ...; and when identifying a particular eigenfunction,
suffix with an~$e$ or~$o$ to indicate $Q$-parity, as in
$\mathcal{B}_e$ or $\mathcal{B}_o$.\footnote{Of course
  the division of eigenfunctions according to dihedral symmetry
  is hardly new. For example, Cureton and Kuttler~\cite{ck1998}
  identify the symmetry classes for the regular hexagon, which can be lined up
  using \mbox{$\mathcal{A}=S0$}, \mbox{$\mathcal{B}_e=S1$},
  \mbox{$\mathcal{C}_o=S2$}, \mbox{$\mathcal{S}'=S3$},
  \mbox{$\mathcal{S}=C0$}, \mbox{$\mathcal{B}_o=C1$},
  \mbox{$\mathcal{C}_e=C2$}, and \mbox{$\mathcal{A}'=C3$}; where the
  \mbox{$S0$--$C3$} are class names in that reference.}

\section*{Some numerical considerations}
When calculating the point-matching determinant, it is important that
all calculations be carried out to relatively high precision. To get
the process started, first calculate a low precision eigenvalue bound
in order to estimate the convergence rate, $\rho$. Then use that to
estimate the \mbox{$N$-value} targeting a high precision bound,
accurate to, say, \mbox{$D$~digits}, i.e., \mbox{$N\leftarrow
  D/\rho$}. To set the precision for this new calculation, use the
rule-of-thumb requiring at least $1.2N$ digits in the intermediate
calculations.  If the process works, the precision was adequate. If
not, try increasing it.

The roots of the determinant must certainly be found to a precision
better than the desired $D$ digits of the eigenvalue, perhaps to
within $1.2D$ digits. I observed that for $D$ beyond a few dozen
digits, the point-matching determinant becomes quite linear in
$\lambda$, so a simple secant-based root-finding algorithm works very
well.

If one pursues hundreds of digits, the computer language must be able
to handle arbitrary-precision calculations efficiently.  Whatever
programming environment you use, be sure that the numerics are good,
especially the fractional-order Bessel functions.\footnote{For
  example, and rather unfortunately, the {\tt maxima} function {\tt
    bessel\_j} does not seem to respect {\tt fpprec}.}

\section*{Examples}
My main project is the regular pentagon, for which I have so-far
calculated the lowest 8139 Dirichlet eigenvalues to at least 60
correctly rounded digits. See Figs.~\ref{fig:stareigenfunctions}
and~\ref{fig:pentagonplots} for several plots, including the highest
one in that set; and---as presented in
Table~\ref{tab:ultraprecise}---some hundred-digit eigenvalue
results. Except for those examples, I find it more interesting to use
other shapes to demonstrate the eigenvalue calculation and bounding
methods.

\subsection*{L-shape}
\begin{figure}[htb!]
\includegraphics[scale=0.6]{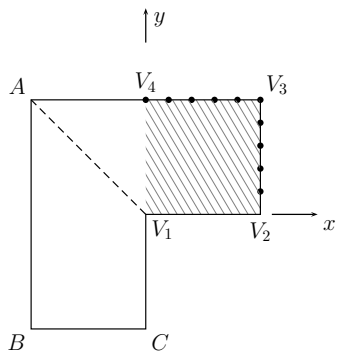}%
\caption{L-shape showing a line of symmetry (even line $V_1A$) and
  region $\Omega$ (shaded square \mbox{$V_1 V_2 V_3 V_4$}, in its
  canonical position) with ten evenly-spaced matching points.}
\label{fig:Lshape}
\end{figure}

This project would not be complete without considering the famous
L-shape formed by joining two unit-edged squares to adjacent edges of
a third, as shown in Fig.~\ref{fig:Lshape}. For this, I shall limit
the scope to calculating the lowest Dirichlet eigenvalue. The first
one-hundred digits of this eigenvalue appear in
Table~\ref{tab:ultraprecise}, and an abbreviated bound is
\begin{equation}
 \lambda_1 = 9.6397238440219410527\cdots 27494610057016061858\,_{4831}^{5300}
\end{equation}
where the leading and trailing digits of the 1001-digit result are
shown. The first 1001 decimal digits of this number are listed as
\url{http://www.oeis.org/A262701}~\cite{oeisorg}.

The 2004 influential ``Reviving the Method of Particular Solutions''
by Betcke and Trefethen~\cite{bt2004}, which, among other things,
examines the L-shape as an example, and apparently provides the best
modern reference for this shape. Despite that, I regress back to 1967
and use FHM as a starting point.

Following FHM, fully exploit symmetry to yield a square region,
$\Omega$, in its canonical position as shown by the shaded region of
Fig.~\ref{fig:Lshape}. This means we can use
Eq.~(\ref{eq:fourierbessel}) with \mbox{$m$-values}
\begin{equation}
              \qquad m_\nu=\frac{2}{3}\left( 2 \left\lfloor\frac{3\nu}{2}\right\rfloor-1\right)
\end{equation}
where $\lfloor x\rfloor$ is $\mathop{\mathrm {floor}}(x)$. The
expression in parentheses gives the required sequence $\{1, 5, 7, 11,
13, ...\}$; i.e., (positive) integers on either side of the integers
in $\{0,6,12,...\}$, which are appropriate for the lowest Dirichlet
eigenmode.

Now, consider a near repeat of the FHM calculation for an
``apples-to-apples'' comparison---this time using arbitrary precision
and much higher \mbox{$N$-values}.  To do this exercise, (a)~choose
equally-spaced matching points on $\partial\Omega_2$ and
$\partial\Omega_3$; (b)~impose
\mbox{$\partial\Psi(k;r,\theta)/\partial\theta=0$} at~$V_2$ and~$V_3$;
and (c)~use~the proper increment of {$\Delta N=2$}.  Then, to
numerically verify the alternating/converging nature, calculate every
approximate eigenvalue for even $N$ from~4 to just
over~260.\footnote{That highest $N$ value yields 100 correct digits in
  the eigenvalue, and altogether, takes about one hour of CPU
  time. Also note that requirement~(b) was imposed by FHM for
  numerical reasons, I impose it here to accurately reproduce their
  results.}  The results for
\mbox{$N=4,6,...,32$,~$254,...,260$} are shown in
Table~\ref{tab:Lshapeloworder}, which illustrates the alternation.

FHM reported the eight-digit $9.6397238\,_{05}^{84}$ (FHM~Table~3)
eigenvalue bound, and they indicate calculations up to $N=26$
(FHM~Table~2).  FHM round-off error appears to become significant for
$N>14$, so use my very low-order \mbox{$N=12$} and~$14$ numbers for
that ``apples-to-apples'' comparison. By inspection, the bound is
$9.6397238\,_{43}^{55}$. This happens to be slightly better, but the
emphasis should be placed on the ease with which this bound is
obtained. Interestingly, FHM data is essentially the same as mine out
to ninth decimal place, so that original FHM (\mbox{$N=12$} and~$14$)
data might be used to bound the eigenvalue to around eight digits
without using Moler-Payne, however, without more precision and higher
N-values, it is not obvious that the alternation may be used to bound
the eigenvalue. Note that the $N$-value used in the FHM Moler-Payne
calculation was not indicated.

\begin{table*}[htb!]
\setlength{\tabcolsep}{1ex}
\caption{Very low-order and a few high-order results for the L-shape
  using the FHM procedure. At least five digits after the last
  matching digit are retained and rounded, and the underlined digits
  guide the eye.  Note the alternation: Very consistently, even and
  odd $N/2$ provide, respectively, upper and lower bounds. The last
  few entries indicate a \mbox{100-digit} calculation (ending in ``\mbox{...788}''), with
the last three showing how to construct the bound. FHM~Table~2
  data are copied here for side-by-side comparison.}
\begin{tabular}{rll}
\toprule
\multicolumn{1}{c}{$N$}& \qquad$\lambda$  & \quad FHM $\lambda^*$ \\
\colrule
4  &  9.\underline{6}58161723        & 9.658161723 \\
6  &  9.63\underline{9}624491        & 9.639624491 \\
8  &  9.6397\underline{2}66319       & 9.639726632 \\
10 &  9.63972\underline{3}70221      & \hphantom{9.63972}3703 \\
12 &  9.639723\underline{8}54826     & \hphantom{9.63972}3855 \\
14 &  9.6397238\underline{4}30369    & \hphantom{9.63972}3844 \\
16 &  9.63972384\underline{4}12442    & \hphantom{9.63972}3844 \\
18 &  9.639723844\underline{0}10281    & \hphantom{9.63972}3845 \\
20 &  9.6397238440\underline{2}33611   & \hphantom{9.63972}3845 \\
22 &  9.63972384402\underline{1}75875  & \hphantom{9.63972}3845 \\
24 &  9.639723844021\underline{9}65466 & \hphantom{9.63972}3844 \\
26 &                    9.639723844021\underline{9}37668 & \hphantom{9.63972}3846 \\
28 & \multicolumn{2}{l}{9.6397238440219\underline{4}15358} \\
30 & \multicolumn{2}{l}{9.6397238440219\underline{4}09820} \\
32 & \multicolumn{2}{l}{9.639723844021941\underline{0}63271}\\
$\vdots$ &&\\
254 & \multicolumn{2}{l}{9.6[$\cdots$]50087967\underline{8}64481\hphantom{6092 }}\\
256 & \multicolumn{2}{l}{9.6[$\cdots$]50087967\underline{8}92906\hphantom{6092 }($\cdots 78\,_{64}^{93}$)}\\
258 & \multicolumn{2}{l}{9.6[$\cdots$]5008796788\underline{7}42560\hphantom{92 }($\cdots 78\,_{87}^{93}$)}\\
260 & \multicolumn{2}{l}{9.6[$\cdots$]5008796788\underline{8}48247\hphantom{92 }($\cdots 788\,_{74}^{85}$)}\\
\botrule
\end{tabular}

\label{tab:Lshapeloworder}
\end{table*}

Except for the P.~Amore et al.~\cite{abfr2015} recent effort, I
believe the best published result is the decade-old, 13-digit
correctly rounded result of Betcke and Trefethen~\cite{bt2004},
\mbox{$9.63972384402\,_{16}^{22}$}. To make a reasonable comparison, they
used the GSVD method, no symmetry reduction, \mbox{$240=4\times 60$}
equal-spaced boundary matching points (excluding all vertices) on four
point-matched segments, 50~randomly-selected interior points, only
15~basis-functions, and machine-precision. That required minimizing
the smallest generalized singular value of a $(50+240)\times 240$
matrix and (to obtain the MP-type bounds) estimating the function
values on the point-matched edge.

Using the current method, which has a better distribution of points
(by including vertices and actually being evenly spaced) and taking
full advantage of symmetry, a similar looking bound 
\mbox{$9.63972384402\,_{17}^{34}$} is achieved by comparing the
\mbox{$N=20$} and~$22$ approximate eigenvalues. This result required a
few seconds of CPU time.  Although not an ``apples-to-apples'' comparison,
it does indicate significantly less numerical effort is required to
obtain a similar result.

Next, I extend the calculations from $N=814$ up to $826$, still using
equally-spaced matching points\footnote{...~but abandoning the
  $\partial\Psi/\partial\theta=0$ requirement at $V_2$ and $V_3$}, to
obtain a \mbox{300-digit} result with an overall convergence rate of
$1/2.75$. This took about a day of CPU time.

Then, I switched to Chebyshev-distributed points (very similar to
those shown in Fig.~\ref{fig:cutsquare} for the cut-square), which
improved the convergence rate to very close to $1/2$. This allowed
calculation of a \mbox{400-digit} result at a similar $N=816$, after
another CPU day.  By extending \mbox{$N$-values} up to about $2100$,
``thousand-digit'' results may be obtained after several weeks of CPU
effort.  This illustrates how exploiting symmetry and judiciously
choosing matching points permits one to extend results to very high
precision.

\subsection*{Cut-square}
\begin{figure}[htb!]
\includegraphics[scale=0.6]{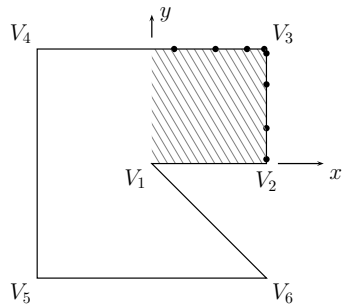}
\caption{Cut-square hexagon and example $N=8$ matching points
  ($N$~even and~$N/2$ on each edge, as shown). There are two matching
  points very close to vertex $V_3$ (one on each edge), and one very
  close to vertex $V_2$. The shaded region is the most convenient and
  smallest $\Omega$ region needed to solve for the non-closed-form
  eigenvalues.}
\label{fig:cutsquare}
\end{figure}

The sole and limited purpose of this example is to illustrate the
bounding method for the apparently difficult shape shown in
Fig.~\ref{fig:cutsquare}, which I call the ``cut square'' since it is
formed by cutting a triangle out of a unit-edged square. It has six
edges and a $7\pi/4$ re-entrant vertex. It was inspired entirely by an
example in Reference~\cite{yh2009}. Those authors, Yuan and He, were
apparently unaware of the then-recent relevant work by Trefethen and
Betcke~\cite{bt2004,tb2006}, but nevertheless provide an interesting
discussion on the L-shape and this cut-square shape---as they attempt
to bound eigenvalues. Of note is that those authors point out that
this shape has no symmetry, but indeed it does.\footnote{Very
  recently, P.~Amore et al.~\cite{abfr2015} also examined this shape
  and discovered other interesting symmetry properties, and were also
  able to calculate very precise Dirichlet eigenvalues. They reported
  40-digit MPS results, which do agree with my results.}

For this polygon, calculations of both Dirichlet and Neumann
eigenvalues are demonstrated. Because of the very high convergence
rate of \mbox{$\rho\approx 1/2$}, with modest effort,
the lowest three eigenvalues for each type are calculated to 200 digits,
with 100-digit truncated values presented in Table~\ref{tab:ultraprecise}.

The first challenge is to identify the symmetries and classify the
eigenmodes. To do this, first consider the Dirichlet modes and expand
the eigenfunction (Eqs.~\ref{eq:fourierbessel} and \ref{eq:basisfunctions}, choosing ``$\sin$'') 
about the re-entrant vertex with
\begin{equation}
            m_\nu=\frac{4\nu}{7} \qquad\mbox{where $\nu=1,2,3,...,N$}.
\end{equation}

By calculating the eigenfunction coefficients (the~$c_\nu$ in
Eq.~(\ref{eq:fourierbessel})), a pattern is revealed where, for a given
eigenfunction, many coefficients are zero.  That lead to a simple
discovery\footnote{This generally identifies a clever technique that
  can be used to empirically discover symmetries if they are not
  obvious or one is not inspired to figure it out from first
  principles.  FHM has indications of this technique.}  
of a symmetry quite analogous to that of the L-shape.
By inspection of the calculated coefficients, it became obvious that
(a)~closed-form modes are selected by choosing integral~$m$, while
(b)~any non-closed-form eigenmode belongs to one of three classes (here labeled A, B, and C)
obtained by choosing positive \mbox{$\nu$-values} on either side of multiples
of seven. To best reveal the pattern, line up the possible \mbox{$m$-values},
\begin{equation}
m =  \textstyle
[0,\!]\,
\overset{A}{\frac{4}{7}},
\overset{B}{\frac{8}{7}},
\overset{C}{\frac{12}{7}},
\overset{C}{\frac{16}{7}},
\overset{B}{\frac{20}{7}},
\overset{A}{\frac{24}{7}},
4, 
\overset{A}{\frac{32}{7}},
\overset{B}{\frac{36}{7}},
\overset{C}{\frac{40}{7}},
\overset{C}{\frac{44}{7}},
\overset{B}{\frac{48}{7}},
\overset{A}{\frac{52}{7}},
8,
\overset{A}{\frac{60}{7}},\cdots
\label{eq:cutsquaremvalues}
\end{equation}
where the zero is shown but used only for the closed-form Neumann solutions.

The fundamental region for the closed-form modes appears to be the
45-90-45 triangle. However, the smallest region in which I can solve
the problem for the non-closed-form classes appears to be twice as large,
$2/7$ of the total area, which can be chosen to be the square region
in the first quadrant. Thus, for the non-closed-form modes, there are
two point-matched edges meeting at right angles---as with the L-shape.

I found that a Chebyshev-like distribution of matching points, equally
divided between the two matching edges on the square in the first
quadrant, and crowded near the 90-degree corners ($V_2$ and $V_3$),
yields excellent convergence rates. Because of these choices, $N$ must
be even, and very specifically,
\begin{equation}
 \mathcal{P}_N:\, \left(x_\mu,y_\mu\right) = \left\{ 
\begin{array}{cl} 
\displaystyle
 \left(\frac12, \frac{1}{4}\left[  1 - \cos \left(\frac{[\mu-\frac12]}{N/2}\,\pi \right)
                      \right] \right) &\qquad\mbox{for}\quad \mu=1,2,\cdots,\frac{N}{2}\\
\displaystyle
 \left(\frac{1}{2}\, \sin \left(\frac{[\mu-\frac12 ]}{N}\,\pi \right)
                       , \frac12 \right)   & \qquad\mbox{for}\quad \mu=\frac{N}{2}+1,\frac{N}{2}+2,\cdots,N
\end{array}\right.
\end{equation}
Incrementing~$\mu$, these matching points start near~$V_2$, and end near
the midpoint of~$V_3$ and~$V_4$, as shown in Fig.~\ref{fig:cutsquare}.

A sweep of the interval from \mbox{$\lambda=0$} to~$100$ reveals the
first three Dirichlet eigenvalues, which belong to class A, B, and C;
respectively. Initial bounds (within~$\pm 1$) of these three are
\mbox{$\lambda_1=36$}, \mbox{$\lambda_2=54$}, and
\mbox{$\lambda_3=74$}; which are then calculated to 200 digits, with
100-digit truncated values presented in Table~\ref{tab:ultraprecise}.
The lowest closed-form Dirichlet
eigenvalue\footnote{$\lambda=(1^2+2^2)(\pi/L)^2=20\pi^2\approx 197.4$,
  where $\Psi(x,y)\propto \sin(\pi x/L)\sin(2\pi y/L)-\sin(2\pi
  x/L)\sin(\pi y/L)$ and $L=1/2$. Not degenerate.}  is relatively high
at $\lambda\approx 197.4$, which appears to be the tenth Dirichlet
eigenvalue.

The discussion for the Neumann modes is nearly identical, except that
one must choose ``$\cos$'' (in Eq.~\ref{eq:basisfunctions}) and, for the
closed-form modes, one must include the $m=0$ term.  Also, this
example is unique among the chosen examples because there are two even
point-matched edges.

Like the Dirichlet results, the lowest three Neumann eigenvalues
belong to symmetry classes~A, B, and~C; respectively. They are
initially bounded (within~$\pm 1$) with \mbox{$\lambda_1=5$},
\mbox{$\lambda_2=12$}, and \mbox{$\lambda_3=18$}; which are then
calculated to 200 digits, with 100-digit truncated values presented in
Table~\ref{tab:ultraprecise}.  The lowest closed-form Neumann
eigenvalue\footnote{$\lambda=(1^2+0^2)(\pi/L)^2=20\pi^2\approx 39.5$,
  where $\Psi(x,y)\propto \cos(\pi x/L)+\cos(\pi y/L)$ and
  $L=1/2$. Not degenerate.}  is $\lambda\approx 39.5$, which is
apparently the next one up.

For these lowest Neumann and Dirichlet eigenvalues, the convergence
rate is quite rapid at about~$1/2$, which means about 400~matching
points are needed to yield the \mbox{200-digit} results.

Prior to this work and that of Ref.~\cite{abfr2015}, the only
published eigenvalue for this shape was the lowest Dirichlet
eigenvalue, and the best bound is the six-digit result,
$\lambda_1=35.6315\,_{15}^{22}$. To get that, Yuan and
He~\cite{yh2009} used a similar MPS starting point, fifty matching
points in total on all but the two adjacent edges forming the ``cut'',
and a far more complicated bound calculation. (Their distribution of
points was unspecified.) To obtain an ``apples-to-apples'' comparison
with that publish result, I temporarily abandoned the symmetry
considerations and repeated the calculation for up to the relatively
low \mbox{$N=49$} matching points (using $N$---a multiple of
seven---Chebyshev matching points, distributed on all four
non-adjacent edges), and was able to write down a nine-digit estimate
\mbox{$\lambda_1=35.6315195\,_{15}^{37}$}, by inspection. That
``low-precision'' result took a few seconds of CPU time.

It is interesting to note that by simply restoring the symmetry
reduction, the convergence rate becomes about three times faster:
At~$50$ expansion terms (i.e., the same numerical effort, requiring a
few CPU seconds), I can write down thirty digits, \mbox{$\lambda_1=
  35.63151951\,7191723095\,2054861420\,_{61}^{79}$}, again, by
inspection.

In hindsight, this problem was not as difficult (to solve) as one
might have been led.  But there are interesting features that do
suggest further research, such as the nature of the symmetry.

\subsection*{Star (regular concave decagon)}
The five-pointed star provides another relatively difficult shape with
which to demonstrate the method. It is also known as a regular concave
decagon or the outline of the pentagram star. The size is such that
points of the star coincide with the vertices of a unit-edged, regular
pentagon.

This shape shares the same symmetry group as the regular
pentagon. Since I observed some regular-pentagon Dirichlet
eigenfunctions develop nodal lines approximating a regular
{\em pentagram\/}---two examples are illustrated in
Figs.~\ref{fig:starandpentagon} and
\ref{fig:DS065a}; I was inspired to examine this shape, at
least for a few low eigenvalues.

A subset of eigenmodes of the regular pentagon do have an approximate
symmetry leading to that observation, however, I reserve that for a
future discussion.  This star can become a research project in and of
itself: It is an interesting and familiar shape for which I have been
unable to find any published eigenvalues for comparison.

\begin{figure}[htb!]
\includegraphics[scale=0.8]{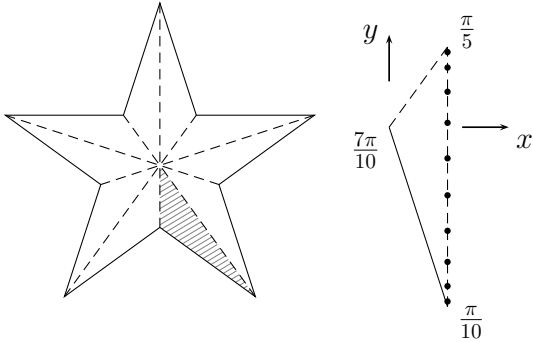}
\caption{(LEFT) Star shape with all even lines connecting opposite
  vertices and (RIGHT) the fundamental region forming the polygon,
  $\Omega$, in its canonical position and with $N=10$ Chebyshev-like
  distributed matching points.}
\label{fig:starB}
\end{figure}

\begin{figure}[htb!]
\includegraphics[scale=1.2]{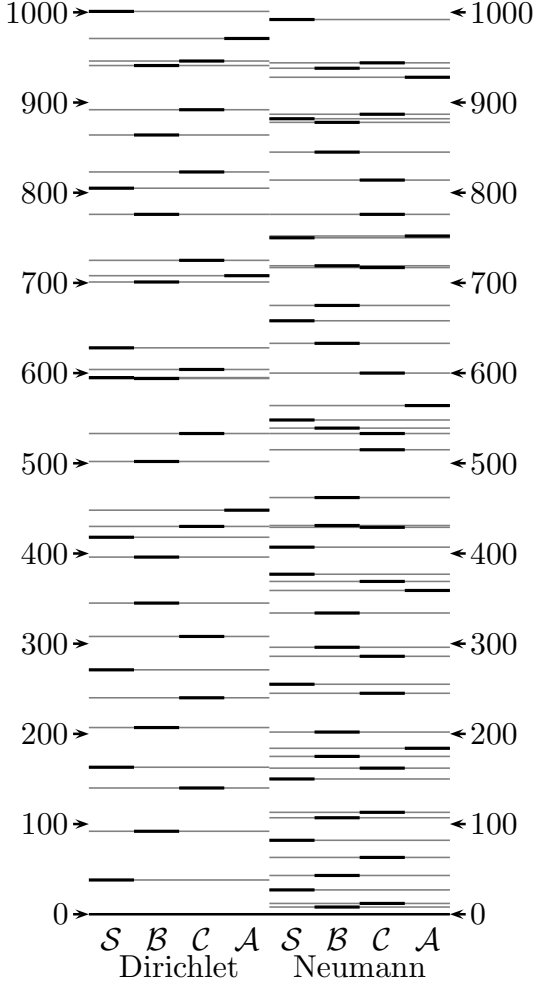}
\caption{The spectrum of the concave decagon (star) for $\lambda < 1002$. This
diagram includes the lowest 79 eigenvalues, which is composed of 31 Dirichlet
and 48 Neumann eigenvalues. They are sorted according to symmetry.
Except that every B and C eigenvalue is doubly degenerate, there are actually no other
degeneracies in this set---even though some come close.}
\label{fig:starSpectrum}
\end{figure}

Like the regular pentagon, all star eigenvalues can sorted into one of
two one-dimensional or two two-dimensional classes.  The
non-degenerate $\mathcal{S}$ and $\mathcal{A}$ eigenvalues require the
triangular fundamental region shown in Fig.~\ref{fig:starB}.  The
doubly-degenerate $\mathcal{B}$ and $\mathcal{C}$ ``mixed'' modes
require the arrow-head shape shown in Fig.~\ref{fig:starA}, which is
formed by reflecting the fundamental triangle through it shortest
edge.

With this problem, unlike some other examples, I could not get the
point-matching method to converge with evenly-spaced matching points,
however, it does work with points crowded near the vertices. A
distribution that yields good convergence rates
is the Chebyshev-like distribution
\begin{equation}
   y_\nu = \frac{1}{2}\left\{\left[ y_{C}+y_{B} \right]
            +\left[ y_{C}-y_{B}\right]\cos\left(\frac{\nu\pi}{n_2+1}\right)\right\}
\end{equation}
for $\nu=1,2,...,n_2$; where $n_2=N$ ($\mathcal{A}$, $\mathcal{S}$) or
$N/2$ ($\mathcal{B}$, $\mathcal{C}$), and where $y_B$ and $y_C$ are
the $y$-values of the matching-edge endpoints.  With this
distribution, there is not a simple alternation of approximate
eigenvalues with incremented values of $N$. Instead, it usually takes
several (about three) steps to alternate. The convergence rate is
nevertheless high enough that the peaks can be interpreted as upper
and lower bounds.

The lowest seventy-nine eigenvalues (including both Dirichlet and
Neumann modes) were calculated for $\epsilon<10^{-20}$, and since so
little is known about this problem, the rich spectrum is simply and
concisely presented in Fig.~\ref{fig:starSpectrum}. Several
representative eigenfunctions are presented in
Fig.~\ref{fig:stareigenfunctions}, and the lowest four Dirichlet and
lowest four Neumann eigenvalues are calculated
to at least 100 digits and presented in Table~\ref{tab:ultraprecise}.

\subsection*{Regular polygons}
Of the regular polygons, it is well known that the equilateral
triangle and the square are fully solved in closed form, and the
regular hexagon has a subset of closed-form solutions (by piecing
together equilateral triangle solutions).  Otherwise, regular polygons
have no closed-form solutions.

Quite surprisingly, (non-closed-form) regular polygon results are
rather fragmentary and spread far and wide in the published
literature. This represents a gaping hole in our collective
understanding. 

To illustrate, consider the regular hexagon. In 1978, Bauer and
Rice~\cite{br1978} calculated the lowest 21 Dirichlet eigenvalues to
at least five or so digits. In 1993, I~\cite{phdthesis} extended that
by a factor of three, to about six digits; and in 1998, Cureton and
Kuttler~\cite{ck1998} almost doubled that tally, and to about eight digits.

The regular pentagon is no better. In 2010, Lanz~\cite{lanz2010}
published what appears to be the most comprehensive list (that I am
aware of) consisting of the ten lowest Dirichlet eigenvalues, accurate
to half-a-dozen digits.

Jumping to the more extreme \mbox{128-sided} regular polygon: In 2004,
Strang and Grinfeld~\cite{sg2004}, and in 2008, Guidotti and
Lambers~\cite{gl2008}, both calculated the ten lowest simple Dirichlet
eigenvalues, accurate to about half-a-dozen digits, the lowest of
which ranges from $5.78319$ to $5.78320$ (adjusting for the area).  In
addition, those authors, and---more
specifically---Oikonomou~\cite{o2010}, offer an asymptotic expansion
in ``$1/\sigma$'' where $\sigma$ is the number of sides on the regular
polygon. Recently, Mark Broady\footnote{Private
  communication.}~\cite{phdthesisBroady} added two more terms to that
expansion, which, for the lowest Dirichlet eigenvalue, is
\begin{equation}
   \frac{\lambda_1}{j_{0,1}^2} = \left\{ 1 + \frac{4\,\zeta(3)}{\sigma^3}
   + \frac{\left[ 12 - 2\, j_{0,1}^2 \right]\,\zeta(5)}{\sigma^5}
   + \frac{\left[ 8 + 4\, j_{0,1}^2 \right]\,\zeta^2(3)}{\sigma^6}
    + \mathcal{O}\left(\frac{1}{\sigma^{7}}\right) \right\}.
\label{eq:regpolasymp}
\end{equation}
where the regular polygon has area equal to~$\pi$,
\mbox{$j_{0,1}\approx 2.4048$} is the first root of the Bessel
function $J_0(x)$, and $\zeta(x)$ is the Riemann-zeta
function.\footnote{Area-$\pi$ regular polygons are preferred since it
  factors out the well-known area dependence on the eigenvalue; and,
  as $\sigma\to\infty$, $\lambda_1\to j_{0,1}^2$.}  Incidentally, and
as far as I know, that has been the only effort to express any
non-closed-form eigenvalues analytically in terms of other ``known''
constants.

To illustrate the method, I find it interesting to calculate the
lowest Dirichlet eigenvalues of $\sigma$-sided regular polygons with
up to $\sigma=130$ at a precision of $\epsilon<10^{-30}$.  The tail
end of that calculation is shown in
Table~\ref{tab:thirtydigitregularpolygon}, with truncated 100-digits
results for the pentagon to the decagon presented in
Table~\ref{tab:ultraprecise}. See also Eq.~(\ref{eq:polygon256}) for a
\mbox{20-digit} result for the 256-sided regular polygon, which was
used to illustrate the notation.

Of note is that above the dodecagon, i.e., $\sigma>12$, the canonical Chebyshev
distribution given by Eq.~(\ref{eq:canonicalchebyshevnodes}) doesn't
work. Instead, using half of it does seem to work for all regular
polygons,
\begin{equation}
   y_\nu = y_{\mathrm max} \sin(\nu\pi/(N+1))
\end{equation}
where $\nu=1,2,...,N$, and where $y_{\mathrm max}$ is the $y$-value of
the highest point of the fundamental triangle in
Fig.~\ref{fig:regularpolygons}\,(CENTER); here $y_{\mathrm
  max}=\cos(\pi/\sigma)$. This distribution crowds matching points
near the ``distant'' acute angle ($\alpha/2=\pi/\sigma$), and spreads
them out near the right angle. It also seems to very reliably yield alternating
approximations at each increment, $\Delta N=1$.

This exercise illustrates specifically that even though exponential
convergence is present, as the number of sides increases, the
convergence rate becomes quite poor.  Beginning at the pentagon,
$\rho\approx 1/1.14$ when $\epsilon\approx 10^{-1000}$. When one
reaches the 130-sided regular polygon, that convergence rate has
dropped to about $\rho\approx 1/28.1$ at the 30-digit result. A rough
estimate for the convergence rate for this problem is $\rho\approx
5/\sigma$, at least for the 30-digit results.

Also of note is that when the convergence rate is good, more precision
in the calculations is needed: The regular pentagon required $1.7N$
digits. Above the decagon, $1.4N$ seemed adequate.

For situations where $\sigma$ (number of polygon sides) is low
(5~or~6), with a little practice, it becomes relatively
straightforward to calculate hundred-digit results, for both Neumann
and Dirichlet boundary conditions, for up to perhaps the lowest ten
thousand eigenvalues. For example, the 8139\textsuperscript{th}
Dirichlet eigenvalue (same symmetry class as the lowest) of the
(unit-edged) regular pentagon to 300 correct digits is bound with
\begin{equation}
\lambda_{8139}=100001.28198274831240\cdots 32644630143129793681\,_{3271}^{5381}
\label{eq:8139}
\end{equation}
which, as indicated, displays the leading and trailing digits.

Turning attention to how far one can extend the precision of these
numbers using a laptop, two ``1000+'' digit eigenvalues are
calculated.  The lowest Dirichlet eigenvalue of the unit-edged regular
pentagon, to 1502 digits,
\begin{equation}
\lambda_1 = 10.996427084559806648 \cdots 82166474404544652968\,^{9936}_{8251}
\end{equation}
was submitted to
\url{http://www.oeis.org/A262823}~\cite{oeisorg}.
This hexagon number required just under one
month of computation time.

Similarly, the lowest Dirichlet eigenvalue of the unit-edged regular
hexagon, to 1001 correct digits,
\begin{equation}
\lambda_1 = 7.1553391339260551282 \cdots 77979918378681828158\,^{1503}_{0889}
\end{equation}
was submitted to
\url{http://www.oeis.org/A263202}~\cite{oeisorg}.
This required about five days of computation time.

\begin{table}[htb!]
  \caption{Thirty-digit results for the lowest Dirichlet eigenvalue
    within a regular polygon with number of sides $\sigma$ and area
    $\pi$.  The last entry is $j_{0,1}^2$, i.e., the corresponding
    eigenvalue of the unit-radius circle. Note how the required number
    of matching points, $\widehat{N}$, increases by about six each
    time to maintain the desired eigenvalue precision. Each of these
    polygon eigenvalues required about a day of CPU.\vspace{0.3ex}}
  \label{tab:thirtydigitregularpolygon}
  \begin{tabular}{clcc}
    \toprule
    $\sigma$ & \qquad $\lambda_1$ ($\epsilon<10^{-30}$) & $\widehat{N}$  & $\epsilon$ \\
    \colrule
    126 &      $5.7831998639169811697955997275\,_{4738}^{5287}$ & 818 & $9.5\times 10^{-31}$\\ 
    127 &      $5.7831995381236804121745520138\,_{6591}^{7147}$ & 824 & $9.6\times 10^{-31}$\\ 
    128 &      $5.78319922243209895698523832013\,_{030}^{555}$  & 831 & $9.7\times 10^{-31}$\\ 
    129 &      $5.78319891645372682901545245421\,_{112}^{643}$  & 837 & $9.2\times 10^{-31}$\\
    130 &      $5.78319861981784749432269771828\,_{013}^{587}$  & 842 & $9.9\times 10^{-31}$\\
    $\infty$ & $5.783185962946784521175995758456$ & & \\
    \botrule
  \end{tabular}
\end{table}

\section*{Acknowledgements}
I wish to thank Alex Barnett for making the specific suggestion to
expand about the non-analytic vertex using fractional-order Bessel
functions, i.e., non-integral $m$-values (private communication,
December 2014). Indeed, after sharing some regular pentagon results
with him, he suggested that instead of expanding about the center of
the regular pentagon, I should expand about one of its vertices. That
one simple, and---in hind-sight---obvious suggestion, immediately
turned my eight-digit calculations into multi-hundred-digit
calculations because of the exponential convergence.

I also wish to thank James Kuttler and Nick Trefethen for suggestions
and encouragement. I am also encouraged by the recent independent efforts
of Mark Broady and Paolo Amore, et al., and wish to thank them for
interesting dialogs.

Of course, this project was made possible by free software, most specifically
{\tt GMP}~\cite{gmp6} and the {\tt GP/PARI}~\cite{PARI2} calculator.

\section*{Conclusion}
By using a method substantially identical to that of FHM, but using
modern hardware and [free] software, and a little patience, I have
demonstrated that it becomes relatively easy to exploit the well-known
``exponential convergence'' to calculate eigenvalues of the Laplacian
to very high precision, very typically hundreds of digits, and about a
thousand digits for the lowest Dirichlet eigenvalue of the L-shape
(1001~digits), regular pentagon (1502~digits), and regular hexagon
(1001~digits).

I believe my unique contribution to this problem includes (a)~the
observation that the sequence of approximate eigenvalues alternates as
one adds matching points---thus yielding an easy and excellent method
to bound eigenvalues, (b)~revealing a very simple way to overcome the numerical
ill-conditioning of the point-matching method and exploit its
advantages to permit easy calculation of those eigenvalues, and
(c)~computing the first hundred-digit and thousand-digit
non-closed-form eigenvalue results for a variety of shapes.

As I explore this interesting and classic problem, I see theoretical
and pedagogical gaps that need to be filled, and a need to compile
results. My hope is to inspire others pick a shape and start
calculating.

\newpage
\section*{Appendix}
\begin{table*}[!htb]
\caption{Selected ultra-precise eigenvalues for the shapes shown in
  Fig.~\ref{fig:allshapes}. Each line has 100 digits in groups of
  ten. Many of these eigenvalues have been calculated to significantly
  higher precision, but all are truncated here to 100 digits. All of
  these reported digits are believed to be correct.}
\newlength{\mygap}
\settowidth{\mygap}{\,}
\newcommand{\M}{\raisebox{-1ex}[0pt][0pt]{\makebox[\mygap][c]{}}}
\begin{tabular}{l}
\toprule
\multicolumn{1}{c}{Lowest Dirichlet eigenvalue of the L-shape (OEIS/A262701)\vspace{0.5ex}}\\
\multicolumn{1}{p{\textwidth}}{%
\mbox{9.639723844\M0219410527\M1145926236\M4823156267\M2895258219\M0645610957\M9700564035\M6478633703\M9072287316\M5008796788}}\\%
\colrule
\multicolumn{1}{c}{Lowest three Dirichlet eigenvalues of the cut-square\vspace{0.5ex}}\\
\mbox{35.63151951\M7191723095\M2054861420\M7765698409\M6719323704\M4187565934\M6885272248\M0828423119\M6945786239\M3580027104}\\
\colrule
\mbox{54.19310844\M4246291974\M1197858564\M7040768914\M7834351054\M6172416365\M9816279545\M9613641230\M6675978946\M4088833202}\\
\colrule
\mbox{73.63330812\M5603834594\M8382867456\M6950026083\M7320383040\M6422202112\M5756253401\M4601480363\M5503338306\M9092168350}\\
\colrule
\multicolumn{1}{c}{Lowest three Neumann eigenvalues of the cut-square\vspace{0.5ex}}\\
\mbox{4.872527692\M6560441293\M9958456263\M8223244356\M0835019173\M6918287660\M7254479092\M9224360057\M0084074146\M2382445186}\\
\colrule
\mbox{11.68901246\M7975646418\M5606630328\M8745685066\M8757229760\M0410850473\M2248922177\M2876660771\M9335431729\M8467748934}\\
\colrule
\mbox{18.41355366\M4057643462\M4365786419\M2874654908\M0932393352\M2677244373\M5428238912\M2726711980\M3455782771\M2026452147}\\
\colrule
\multicolumn{1}{c}{Lowest four Dirichlet eigenvalues of the star\vspace{0.5ex}}\\
\multicolumn{1}{p{\textwidth}}{%
\mbox{38.16467784\M9021956120\M7064401254\M4909302761\M7878759405\M6032173115\M9893430940\M1381888098\M1660705230\M8934649411}}\\
\colrule
\mbox{91.52297660\M0087650110\M1874055059\M4495578448\M7274004759\M3637987954\M9616910512\M9744672654\M0219433040\M9841939843}\\
\colrule
\mbox{140.1751508\M7431829715\M7695648583\M4562984886\M8950461859\M0551648700\M0199831865\M5965928820\M7377894288\M1434147089}\\
\colrule
\multicolumn{1}{p{\textwidth}}{%
\mbox{163.1376917\M3762835226\M6351068194\M6753674371\M1941271352\M4052991614\M0611124064\M5716586029\M2880248852\M0079395340}}\\
\colrule
\multicolumn{1}{c}{Lowest four Neumann eigenvalues of the star\vspace{0.5ex}}\\
\mbox{8.142790964\M1219464723\M3479298481\M4692952302\M4540073157\M4281385071\M7750378042\M0287187929\M1379273329\M1003542992}\\
\colrule
\mbox{12.39876832\M7444927056\M0165051825\M7451412024\M6307375894\M2540568120\M9312929258\M9363195125\M1866820544\M1275209247}\\
\colrule
\mbox{27.35593497\M9775790508\M6023809638\M0191310466\M0644900698\M7296023165\M9815796111\M1183248158\M1641731648\M2599322915}\\
\colrule
\mbox{43.45373653\M4633280721\M7484231722\M7400872562\M0705962339\M9821299710\M0656883686\M6791890220\M6208479355\M9294072629}\\ 
\colrule
\colrule
\multicolumn{1}{c}{Lowest Dirichlet eigenvalue of regular polygons (area=$\pi$), sides $\sigma=5,6,7,8,9,10$\vspace{0.5ex}}\\
\mbox{6.022137932\M0426338782\M9800871005\M4242967005\M3053404485\M5798207884\M6736673784\M8013486862\M2358041642\M6117672672}\\
\colrule
\mbox{5.917417831\M6136612156\M8857457683\M8961545008\M2860040929\M3411904080\M5009620004\M4903613632\M2387404110\M3755146850}\\
\colrule
\mbox{5.866449312\M6559858577\M1247494175\M8841084242\M7349136980\M5378445601\M2986066686\M0112738387\M6297693674\M2836505116}\\
\colrule
\mbox{5.838491433\M5924428505\M1664037956\M3815784836\M7571520259\M6209417605\M2525113941\M2694626037\M2422792179\M6589151583}\\
\colrule
\mbox{5.821826802\M2702657317\M3554644371\M6945921671\M7867646205\M8910792511\M5921840916\M6829763226\M0595285520\M2545987984}\\
\colrule
\mbox{5.811260359\M2191160227\M8881646881\M1164623442\M1581749002\M3191327488\M5882804261\M2940209711\M4160422652\M9269462629}\\
\botrule
\end{tabular}

\vspace{1in}
\label{tab:ultraprecise}
\end{table*}
\bibliographystyle{plain}
\bibliography{references}
\end{document}